\documentclass[a4,psfig,english,11pt,epsf,portrait]{article}
\usepackage{amsmath,amsfonts,latexsym,amscd,amssymb,theorem}
\usepackage{epsfig}
\usepackage{amsmath}
\usepackage{color}
\definecolor{purple}{rgb}{0.65, 0, 1}
\definecolor{orange}{rgb}{1,.5,0}
\definecolor{brown}{rgb}{.9,.73,.26}

\def\R{\hbox{\bf R}}
\def\Z{\hbox{\bf Z}}
\def\S{{\mathbb S}}

\def\N{\hbox{\bf N}}

\def\b{\beta}

\def\O{\Omega}

\def \proof{{\sl {\bf Proof.}}}

\def\<{\langle}
\def\>{\rangle}
\newcommand{\ba}{\begin{eqnarray}}
\newcommand{\ea}{\end{eqnarray}}

\newtheorem{thm}{Theorem}[section]

\newtheorem{lem}[thm]{Lemma}

\newtheorem{theorem}[thm]{Theorem}
\newtheorem{definition}[thm]{Definition}
\newtheorem{lemma}[thm]{Lemma}
\newtheorem{proposition}[thm]{Proposition}
\newtheorem{corollary}[thm]{Corollary}
\newtheorem{rem}[thm]{Remark}

\newcommand{\eps}{\varepsilon}
\numberwithin{equation}{section}

\renewcommand{\R}{{\mathbb R}}
\renewcommand{\Z}{{\mathbb Z}}
\renewcommand{\N}{{\mathbb N}}

\begin{document}
\title{\bf The Peierls-Nabarro model as a limit of a Frenkel-Kontorova model}
\author{
\normalsize\textsc{ A. Z. Fino\footnote{D\'{e}partement de
Math\'{e}matiques,
    Laboratoire de math\'ematiques appliqu\'ees, UMR CNRS 5142, Universit\'e de Pau et des Pays
de l'Adour, 64000 Pau, France. E-mail:
    ahmad.fino01@gmail.com
\newline \indent $\,\,{}^{a}$LaMA-Liban,
Lebanese University, P.O. Box 37 Tripoli, Lebanon.
\newline \indent $\,\,{}^{**}$Invited professor at LaMA-Liban, Lebanese University, P.O. Box 37 Tripoli, Lebanon.
\newline \indent $\,\,{}^{b}$Lebanese University, Faculty of Sciences, Mathematics Department, Hadeth, Beirut,
Lebanon. E-mail: ibrahim@cermics.enpc.fr
\newline \indent $\,\,{}^{c}$Universit\'{e} Paris-Est, CERMICS,
Ecole des Ponts, 6 et 8 avenue Blaise Pascal, Cit\'e Descartes
Champs-sur-Marne, 77455 Marne-la-Vall\'ee Cedex 2, France. E-mail:
monneau@cermics.enpc.fr
 \newline \indent $\,\,{}^{d}$ D\'epartement de Math\'ematiques, Universit\'e de La
Rochelle,  17042 La Rochelle, France. }$^{\,\,\, ,\,a,\,d}$, H.
Ibrahim$^{a,\,b}$, R. Monneau$^{**,\,c}$
  }}

\vspace{20pt}

\maketitle


\centerline{\small{\bf{Abstract}}} \noindent{\small{We study a
generalization of the fully overdamped Frenkel-Kontorova model in
dimension {\color{black}$n\geq 1.$} This model describes the
evolution of the position of each atom in a crystal, and is
mathematically given by an infinite system of coupled first order
ODEs. We prove that for a suitable rescaling of this model, the
solution converges to the solution of a Peierls-Nabarro model, which
is a coupled system of two PDEs (typically an elliptic PDE in a
domain with an evolution PDE on the boundary of the domain). This
passage from the discrete model to a continuous model is done in the
framework of viscosity solutions.}} \hfill\break

\noindent {\it \small \bf MSC:} $49$L$25;$ $30$E$25$

\section{Introduction}\label{sec1}
In this paper we are interested in two models describing the
evolution of defects in crystals, called dislocations. These two
models are the Frenkel-Kontorova model and the Peierls-Nabarro model. The
Frenkel-Kontorova model is a discrete model which describes the
evolution of the position of atoms in a crystal. On the contrary,
the Peierls-Nabarro model is a continuous model where the
dislocation is seen as a phase transition. The main goal of the
paper is to show rigorously how the Peierls-Nabarro model can be
obtained as a limit of the Frenkel-Kontorova model after a suitable
rescaling.
\subsection{Peierls-Nabarro model}
Let us start to present the Peierls-Nabarro model. We set
$$\Omega = \{x=(x_{1},..., x_{n})\in \R^{n},\; x_{n}>0\},$$
and for a time $0<T\leq+\infty,$ we look for solutions $u^0$ of the
following system with $\beta\geq0$:
\begin{equation}\label{PN}
\left\{
\begin{aligned}
& \beta u^0_t(x,t)=\Delta u^0(x,t) &(x,t)\in&\,\Omega\times (0,T)&\\
& u^0_t(x,t)=F(u^{0}(x,t))+\frac{\partial u^0}{\partial
    x_{n}}(x,t)  &(x,t)\in&\,\partial{\Omega}\times (0,T),&\\
 \end{aligned}
\right.
\end{equation}
where the boundary $\partial\O$ is defined by:
$$\partial\Omega = \{x=(x_{1},..., x_{n})\in \R^{n},\; x_{n}=0\},$$
and the unknown $u^{0}(x,t)\in \R$ is a scalar-valued
function with the initial data
\begin{equation}\label{Init_cond_PN}
u^{0}(x,0)=u_{0}(x)\quad \mbox{for} \quad
\left\{
\begin{aligned}
& x\in \overline{\O}:=\partial\O\cup\O\quad &\mbox{if}& \quad \beta>0\\
& x\in \partial\O\quad
&\mbox{if}& \quad \beta= 0.
\end{aligned}
\right.
\end{equation}
We assume the following conditions on the function
$F:\;\R\rightarrow\R$ and on the initial data $u_0:$
\begin{equation}\label{Reg1}
F\in {\color{black}W^{2,\infty}(\R)}\quad\mbox{and}\quad \left\{
\begin{aligned}
& u_0\in W^{3,\infty}(\overline{\O})\quad &\mbox{if}& \quad \beta>0\\
& u_0\in W^{3,\infty}(\partial\O)\quad
&\mbox{if}& \quad \beta= 0.
\end{aligned}
\right.\end{equation} For simplicity we have taken a high regularity
on the initial data, but this condition can be weaken. The classical
Peierls-Nabarro model corresponds to the case $\beta=0$ with a
$1$-periodic function $F$ (see Section~\ref{sec2}). For a derivation
and study of the Peierls-Nabarro model, we refer the reader to
\cite{hl} (and \cite{p,n} for the original papers and \cite{N} for a
recent review by Nabarro on the Peierls-Nabarro model). Let us
mention that a physical and numerical study of the evolution problem
(\ref{PN}) has been treated in \cite{MBW1998}.

In this paper we consider the general case $\beta\geq0,$ since the
case $\beta=0$ is natural in the Frenkel-Kontorova model (see
Section \ref{sec2}), and mathematically the case $\beta>0$ is easier
to study. In the special case $\beta=0$, this problem can be
reformulated (at least at a formal level) as a nonlocal evolution
equation written on $\partial \Omega$ in any dimension. See in
particular \cite{GM} for such a reformulation in dimension $n=2$ and
a limit of the Peierls-Nabarro model to the discrete dislocation
dynamics (see also \cite{MP} for a homogenization result of the
Peierls-Nabarro model). We refer the reader to
\cite{Cabre-SolaMorales} for the study of stationary solutions of
our model when $F=-W'$.

\subsection{Frenkel-Kontorova model}
We now present the Frenkel-Kontorova model. This is a discrete model
which contains a small scale $\varepsilon>0$ that can be seen as the
order of the distance between atoms. We set the discrete analogue of
$\O$:
$$\O^{\eps} = (\eps \Z)^{n-1} \times \eps(\N\setminus \{0\})$$
with discrete boundary
$$\partial \O^{\eps} = (\eps \Z)^{n-1} \times \{0\},$$
and we set
$$\overline{\O}^{\eps}= \O^{\eps} \cup\partial \O^{\eps}.$$
We look for solutions $u^\varepsilon(x,t)$ with
$x\in\overline{\O}^\varepsilon,$ $t\geq0,$ of the following system for $\beta \ge 0$:
\begin{equation}\label{FK}
\left\{
\begin{aligned}
&  \beta
u^\varepsilon_t(x,t)=\Delta^\varepsilon [u^\varepsilon](x,t),
   &(x,t)\in&\,\Omega^{\varepsilon}\times (0,T)\\
&
u^\varepsilon_t(x,t)=F(u^{\eps}(x,t))+D^\varepsilon[u^\varepsilon](x,t),
&(x,t)\in&\,\partial{\Omega^{\varepsilon}}\times (0,T),
\end{aligned}
\right.
\end{equation}
with initial data
\begin{equation}\label{Init_cond_FK}
u^{\eps}(x,0)=u_{0}(x)\quad \mbox{for} \quad
\left\{
\begin{aligned}
& x\in \overline{\O}^{\eps}\quad
&\mbox{if}& \quad \b> 0\\
& x\in \partial \O^{\eps} \quad &\mbox{if}& \quad \beta=0,
\end{aligned}
\right.
\end{equation}
where we have used the following notation for the discrete analogue
of $\Delta$ and $\frac{\partial}{\partial x_n}$ respectively:
\begin{equation}\label{deltaeps}
\left\{
\begin{aligned}
& \Delta^\varepsilon [u^\varepsilon](x,t)=\frac{1}{\varepsilon^2}
\sum_{{\color{black}y}\in\mathbb{Z}^{n},\,|{\color{black}y}|=1} (u^\varepsilon(x+\varepsilon {\color{black}y},t)-u^\varepsilon(x,t))\\
&
D^\varepsilon[u^\varepsilon](x,t)=\frac{1}{\eps}\sum_{{\color{black}y}\in\mathbb{Z}^n,\,|{\color{black}y}|=1,\,{\color{black}y}_{n}\geq0}(u^\varepsilon(x+\varepsilon
{\color{black}y},t)-u^\varepsilon(x,t)).
\end{aligned}
\right.
\end{equation}
In this model, the function  $u^\varepsilon(x,t)$ describes the
scalar analogue of the position of the atom of index $x$ in the
crystal. The classical fully overdamped Frenkel-Kontorova model is a
one dimensional model. It corresponds to the particular case $n=2$
where $u^\varepsilon(x,t)=0$ if $x\in\O^\varepsilon,$ and therefore
only the $u^\varepsilon(\varepsilon  {\color{black}x},t)$ for
${\color{black}x}\in\Z\times\{0\}$ can be non trivial (we refer the
reader to \cite{BK}). Remark that in the case $\beta>0,$ we can
consider classical solutions using the Cauchy-Lipschitz theory,
while at least in the case $\beta=0,$ it is convenient to deal with
viscosity solutions (see Section \ref{sec3}).

More generally, Frenkel-Kontorova are also used for the description
of vacancy defects at equilibruim, see \cite{HMR2007}. See also
\cite{HMR2006}, where the authors study the problem involving a
dislocation inside the interphase between two identical lattices.
Their model corresponds to our model (\ref{FK}) at the equilibrium
with $F=-{W}'$ where the potential $W$ is a cosine function. For
other 2D FK models, see \cite{CB,CCHO}. For homogenization results
of some FK models, we refer the reader to \cite{FIM2}.

\subsection{Main result}

Our main result is the following theorem which establishes
rigorously (for the first time up to our knowledge) the link between
the two famous physical models: the Frenkel-Kontorova model and the
Peierls-Nabarro model.
\begin{theorem}$(\textbf{Existence, uniqueness and convergence})$\label{theo1}\\
Let $\eps>0,$ $\beta\geq0$ and $0<T\leq+\infty.$ Under the condition
$(\ref{Reg1})$ there exists a unique discrete viscosity solution
$u^{\eps}$ of $(\ref{FK})$-$(\ref{Init_cond_FK}).$ Moreover, as
$\eps\rightarrow 0$, the sequence $u^{\eps}$ converges to the unique
bounded viscosity solution $u^{0}$ of $(\ref{PN})$-$(\ref{Init_cond_PN}).$ The
convergence $u^{\eps}\rightarrow u^{0}$ has to be understood in the
following sense: for any compact set $K\subset
\overline{\Omega}\times [0,T)$, we have:
$$\|u^{\eps}-u^{0}\|_{L^{\infty}(K \cap (\overline{\Omega}^\eps\times
  [0,T)))} \rightarrow 0\quad\mbox{as}\quad \eps \rightarrow 0.$$
\end{theorem}
A formal version of this result has been announced in \cite{ELIM}.
We refer the reader to Section $\ref{sec3}$ for the definition of
viscosity solutions. The proof of convergence is done in the
framework of viscosity solutions, using the half-relaxed limits. The
uniqueness of the limit $u^0$ follows from a comparison principle
that we prove for $(\ref{PN})$, using in particular a special test
function introduced by Barles in \cite{Barles}. {\color{black} Most
of the difficulties arise here from the unboundedness of the domain,
and in the non standard case $\beta =0$, where the evolution
equation is on the boundary of the domain instead of being in the
interior of the domain. } Let us notice that because system
(\ref{FK}) can be seen as a discretization scheme of system
(\ref{PN}), then Theorem \ref{theo1} can be interpreted as a
convergence result for such a scheme. In the literature on numerical
analysis of finite difference schemes, we find other results of
convergence for different equations (see for instance \cite{BS} and
\cite{BJ}). Let us mention that an estimate on the rate of
convergence of $u^\varepsilon$ to $u^0$ is still an open question
for our system.

{\color{black}
\begin{rem}
From the proof of Theorem \ref{theo1}, it is easy to see that the
result still holds true if $F$ is replaced by a sequence
$(F^\varepsilon)_\varepsilon$ of functions which converges in
$W^{2,\infty}(\R)$. For such an example of application, see
(\ref{eq::fea}) at the end of Section \ref{sec2}; see also
\cite{ELIM}.
\end{rem}
}

\subsection{Organization of the paper}

This paper is organized as follows. In Section~\ref{sec2}, we
present the physical motivation to our problem. In
Section~\ref{sec3}, we present the definitions of viscosity
solutions for the discrete and continuous problem
$(\ref{FK})$-$(\ref{Init_cond_FK})$ and
$(\ref{PN})$-$(\ref{Init_cond_PN})$  respectively.
Section $\ref{sec6}$ is dedicated to construct uniform barriers of the
solution $u^\eps$ of $(\ref{FK})$-$(\ref{Init_cond_FK})$.
Using those barriers, we prove in Section~\ref{sec5} the existence of a solution for the
discrete problem $(\ref{FK})$-$(\ref{Init_cond_FK})$.
Section~\ref{sec4} is devoted to the proof of Theorem~\ref{theo1}.
Sections~\ref{sec7} and \ref{sec8} are respectively devoted to the proofs of the
comparison principle for the discrete problem
$(\ref{FK})$-$(\ref{Init_cond_FK})$ and the continuous problem
$(\ref{PN})$-$(\ref{Init_cond_PN})$, that were presented in Section \ref{sec3} and used in the whole paper.
Finally, in  the Appendix, we give a convenient corollary of Ishii's lemma, that we use in Section \ref{sec8}.


\section{Physical motivation}\label{sec2}

\subsection{Geometrical description}

We start  to call $(e_1,e_2,e_3)$ an orthonormal basis of the three
dimensional space. In the corresponding coordinates, we consider a
three-dimensional crystal where each atom is initially at the
position of a node ${\color{black}I}$ of the lattice
$$\Lambda=\Z^2\times \left(\frac{1}{2} +\Z\right).$$
This lattice is simply obtained by a translation along the vector
$\frac12 e_3$ of the lattice $\Z^3$, and will be more convenient for
the derivation of our model. We assume that each atom
${\color{black}I}\in\Lambda$ of the crystal has the freedom to move
to another position ${\color{black}I}+ u_{{\color{black}I}}e_2$
where $u_{{\color{black}I}}\in \R$ is a unidimensional displacement.
In particular, we will be able to describe dislocations only with
Burgers vectors that are multiples of the vector $e_2$ (see
\cite{hl} for an introduction to dislocations and a definition of
the Burgers vector). Moreover we introduce the general notation
$$\bar{{\color{black}I}}=({\color{black}I}',-{\color{black}I}_3) \quad \mbox{for all}\quad {\color{black}I}=({\color{black}I}',{\color{black}I}_3)\in \Z^2\times \left(\frac12 +\Z\right) \quad \mbox{with}\quad
{\color{black}I}'=({\color{black}I}_1,{\color{black}I}_2)\in\Z^2$$
and will only consider antisymmetric displacements
$u_{{\color{black}I}}$, i.e. satisfying
\begin{equation}\label{eq::anti}
u_{\bar{{\color{black}I}}}=-u_{{\color{black}I}} \quad \mbox{for
all}\quad {\color{black}I}\in\Lambda.
\end{equation}
We also introduce the discrete gradient
$$(\nabla^d u)_{{\color{black}I}}=\left(\begin{array}{l}
u_{{\color{black}I}+e_1}-u_{{\color{black}I}}\\
u_{{\color{black}I}+e_2}-u_{{\color{black}I}}\\
u_{{\color{black}I}+e_3}-u_{{\color{black}I}}
\end{array}\right).$$
Generally, the core of a dislocation is localized where the discrete
gradient is not small. In our model, screw and edge dislocations
will be represented as follows.

\noindent {\bf Screw dislocation}\\
For
${\color{black}I}=({\color{black}I}_1,{\color{black}I}_2,{\color{black}I}_3)$,
we can consider a dislocation line parallel to the vector $e_2$,
with a displacement $u_{{\color{black}I}}$ independent on
${\color{black}I}_2$, i.e. such that
$u_{{\color{black}I}}=s({\color{black}I}_1,{\color{black}I}_3)$. We
will assume for instance that
$$\left\{\begin{array}{l}
s({\color{black}I}_1,+\frac12)=-s({\color{black}I}_1,-\frac12) \to 0
\quad \mbox{as}\quad {\color{black}I}_1 \to
-\infty\\
\\
s({\color{black}I}_1,+\frac12)=-s({\color{black}I}_1,-\frac12) \to
\frac12  \quad \mbox{as}\quad {\color{black}I}_1 \to +\infty.
\end{array}\right.$$
Moreover, if there is no applied stress on the crystal, then it is
reasonable to assume that
$$\mbox{dist }\left((\nabla^d u)_{{\color{black}I}}, \Z^3 \right) \to 0 \quad \mbox{as}\quad
|({\color{black}I}_1,{\color{black}I}_3)|\to +\infty$$ which means
that the crystal is perfect far from the core of the dislocation.

\noindent {\bf Edge dislocation}\\
For
${\color{black}I}=({\color{black}I}_1,{\color{black}I}_2,{\color{black}I}_3)$,
we can consider a dislocation line parallel to the vector $e_1$,
with a displacement $u_{{\color{black}I}}$ independent on
${\color{black}I}_1$, i.e.
$u_{{\color{black}I}}=e({\color{black}I}_2,{\color{black}I}_3)$. We
will assume for instance that
$$\left\{\begin{array}{l}
e({\color{black}I}_2,+\frac12)=-e({\color{black}I}_2,-\frac12) \to 0
\quad \mbox{as}\quad {\color{black}I}_2 \to
-\infty\\
\\
e({\color{black}I}_2,+\frac12)=-e({\color{black}I}_2,-\frac12) \to
\frac12  \quad \mbox{as}\quad {\color{black}I}_2 \to +\infty.
\end{array}\right.$$
Moreover, if there is no applied stress on the crystal, then it is
reasonable to assume that
$$\mbox{dist}\left((\nabla^d u)_{{\color{black}I}}, \Z^3 \right) \to 0 \quad \mbox{as}\quad
|({\color{black}I}_2,{\color{black}I}_3)|\to +\infty.$$ Remark that
this model can also describe more generally curved dislocations,
which are neither screw, nor edge, but are mixed dislocations.

\subsection{Energy of the crystal}

We assume that each atom $I$ is related to its nearest neighbors
${\color{black}J}$ by a nonlinear spring, whose force is derived
from a smooth potential {\color{black}
$W_{{\color{black}I}{\color{black}J}}$. Then formally the full
energy of the crystal after a displacement
$u=(u_{{\color{black}I}})_{{\color{black}I}\in\Lambda}$ is
$$E(u)=\frac12  \sum_{{\color{black}I},{\color{black}J}\in \Lambda,\ |{\color{black}I}-{\color{black}J}|=1} W_{{\color{black}I}{\color{black}J}}(u_{{\color{black}I}}-u_{{\color{black}J}}).$$
We assume that the dislocation cores are only included in the double
plane ${\color{black}I}_3=\pm \frac12$. We also assume that
$$W_{{\color{black}I}{\color{black}J}}(a)=\left\{\begin{array}{ll}
\varepsilon W(a) & \quad \mbox{if}  \quad|{\color{black}I}_3|=|{\color{black}J}_3|=\frac12\\
W(a) & \quad \mbox{if}  \quad|{\color{black}I}_3|\not=\frac12 \quad
\mbox{or}\quad  |{\color{black}J}_3|\not=\frac12.
\end{array}\right.$$
Here $\varepsilon>0$ is a small parameter. This means that the
springs lying between the two planes ${\color{black}I}_3=\pm
\frac12$ are very weak in comparison to the other springs. This will
allow the core of the dislocation to spread out on the lattice as
$\varepsilon\to 0$ and will allow us to recover the Peierls-Nabarro
model in this limit, after a suitable rescaling.

Notice that, to be compatible with the antisymmetry (\ref{eq::anti})
of the crystal, it is reasonable to assume that
$$W(-a)=W(a) \quad \mbox{for all}\quad a\in \R.$$
Recall that each line of atoms ${\color{black}I}^0+\Z e_2$ only
contains atoms of identical nature (i.e. surrounded by a similar
configuration of springs). Therefore for $u_{{\color{black}I}} =
s({\color{black}I}_1,{\color{black}I}_3)$, a lattice of atoms at the
position ${\color{black}I}+u_{{\color{black}I}}e_2$ and another
lattice of atoms at the position
${\color{black}I}+u_{{\color{black}I}}e_2 +
k_{({\color{black}I}_1,{\color{black}I}_3)}e_2$ with arbitrary
$k_{({\color{black}I}_1,{\color{black}I}_3)} \in\Z$ are completely
equivalent. } Therefore their energy should be the same, and it is
then natural to assume that the potential $W$ is $1$-periodic, i.e.
satisfies
$$W(a+1)=W(a) \quad \mbox{for all}\quad a\in \R.$$
We will also assume that non-deformed lattices minimize globally their
energy, i.e. satisfy
$$W(\Z)=0 \quad \quad \mbox{and}\quad W>0 \quad \mbox{on}\quad \R\backslash
\Z.$$ Let us now assume that a constant global shear stress is
applied on the crystal such that it creates a global shear of the
crystal in the direction $e_2$ with respect to the coordinate
${\color{black}I}_3$. This means that there exists a constant
$\tau\in\R$ such that
$$\mbox{dist }\left((\nabla^d u)_{{\color{black}I}} -\tau e_3, \Z^3\right) \to 0$$
as ${\color{black}I}$ goes far away from the core of the
dislocation. Remark that this assumption is compatible with the
antisymmetry (\ref{eq::anti}) of $u$. {\color{black}If there is no
dislocations, we can simply take for instance
$u_{{\color{black}I}}=\tau {\color{black}I}_3$.}

\subsection{Fully overdamped dynamics of the crystal}

{\color{black} The natural dynamics should be given by Newton's law
satisfied by each atom. This dynamics is very rich and for certain
shear stress $\tau$, it is known (see \cite{CF,HW}) in 2D lattices
that certain edge dislocations can propagate with constant mean
velocity. This is due to the fact that part of the energy is lost by
radiation of sound waves in the crystal. This phenomenon is similar
to the effective drag force created by the surrounding fluid on a
boat or an air plane. The resulting behaviour is a kind of
dissipative dynamics that we modelise here by a fully overdamped
dynamics of the crystal. See \cite{KT1,KT2,KT}) for a fundamental
justification of overdamped type dynamics based on explicit
computations in a 1D Hamiltonian model. For general physical
justifications of the dissipative effects in the motion of
dislocations, see also \cite{hl,AI}. }

We recall that we assume that the dislocation cores are only
contained in the double plane ${\color{black}I}_3=\pm \frac12$. For
this reason, we will artificially distinguish the dynamics inside
this double plane and outside this double plane. We consider the
following fully overdamped dynamics (where the velocity of each atom
is proportional to the force deriving from the energy) which is
written formally as:
\begin{equation}\label{eq::dyna}
\alpha_{{\color{black}I}} \dot{u}_{{\color{black}I}}= -
\nabla_{u_{{\color{black}I}}} E(u) \quad \mbox{with}\quad
\alpha_{{\color{black}I}}= \left\{\begin{array}{ll}
1 & \quad \mbox{if}\quad |{\color{black}I}_3|= \frac12\\
\alpha  & \quad \mbox{if}\quad |{\color{black}I}_3|\not=  \frac12,
\end{array}\right.
\end{equation}
for some constant $\alpha \ge 0$. Here $\dot{u}_{{\color{black}I}}$
denotes the time derivative of the displacement
$u_{{\color{black}I}}$ and
$$\nabla_{u_{{\color{black}I}}} E(u) := \sum_{{\color{black}J}\in\Lambda,\ |{\color{black}J}-{\color{black}I}|=1} W'_{{\color{black}I}{\color{black}J}}(u_{{\color{black}I}}-u_{{\color{black}J}}).$$
Remark in particular that this dynamics preserves the antisymmetry
(\ref{eq::anti}) of $u$. {\color{black} For slow motion of
dislocations (i.e. small velocity with respect to the velocity of
sound in the crystal), it is reasonable to assume that $\alpha=0$,
i.e. the lattice is instantaneously at the equilibrium outside the
double plane ${\color{black}I}_3=\pm \frac12$.

We also assume that the potential  is harmonic close to its minima, i.e. satisfies
$$W(a)=\frac{a^2}{2} \quad \mbox{for}\quad |a|\le \delta < \frac12$$
and assume that the strain is small enough in the crystal outside
the double plane $I_3=\pm \frac12$, i.e. we assume that
$$|u_{{\color{black}J}}-u_{{\color{black}I}}|\le \delta \quad \mbox{for any}\quad {{\color{black}J}},{\color{black}I}\in \Lambda \quad
\mbox{such that}\quad |{\color{black}J}-{\color{black}I}|=1 \quad
\mbox{and}\quad |I_3| \not= \frac12 \quad \mbox{or}\quad
|J_3|\not=\frac12.$$ This assumption allows us, in the region
${\color{black}I}_3\not=\pm \frac12$, to consider forces that can be
expressed linearly in terms of the displacement. We also set
$$\tau=\varepsilon \sigma \quad \mbox{for some}\quad \sigma\in
\R \quad \mbox{and}\quad
v_{{\color{black}I}}=u_{{\color{black}I}}-\varepsilon \sigma
{\color{black}I}_{3}.$$ Using the antisymmetry of the solution
(\ref{eq::anti}), we deduce that we can rewrite the dynamics
(\ref{eq::dyna}) with $\alpha=0$ as
$$\left\{\begin{array}{ll}
\displaystyle{0= \sum_{{\color{black}J}\in\Lambda,\
|{\color{black}J}-{\color{black}I}|=1}
  (v_{{\color{black}J}}-v_{{\color{black}I}})} & \displaystyle{\quad \mbox{for}\quad {\color{black}I}_3 > \frac12}\\
&\\
\displaystyle{\dot{v}_{{\color{black}I}} = \varepsilon\sigma
-\varepsilon W'(2v_{{\color{black}I}}+\varepsilon
  \sigma)+
  \sum_{{\color{black}J}\in\Lambda,\ |{\color{black}J}-{\color{black}I}|=1,\ {\color{black}J}_3\ge \frac12} (v_{{\color{black}J}}-v_{{\color{black}I}})}
&\displaystyle{\quad \mbox{for}\quad {\color{black}I}_3 = \frac12}.
\end{array}\right.$$
Then we see that
$$u^\varepsilon(x,t)=v_{(\frac{x}{\varepsilon}+\frac12)}\left(\frac{t}{\varepsilon}\right)$$
solves (\ref{FK}) with $\beta=0$ and $F$ replaced by
\begin{equation}\label{eq::fea}
F^\varepsilon(a)=\sigma -W'(2a +\varepsilon \sigma).
\end{equation}

\begin{rem}
At the level of modeling, we could also consider more general lattices than $\Z^n$
with more general nearest neighbors interactions,
but this case is not covered by the result of Theorem \ref{theo1} and would require a specific work.
\end{rem}
}

\section{Viscosity solutions }\label{sec3}

In this section we present the notion of viscosity solutions and
some of their properties for the discrete problem
$(\ref{FK})$-$(\ref{Init_cond_FK})$  and then for the continuous
problem $(\ref{PN})$-$(\ref{Init_cond_PN}).$ For the classical
notion of viscosity solutions, we refer the reader to \cite{BC-D},
\cite{Bar}, \cite{CIL}.

\subsection{Viscosity solutions for the discrete problem}\label{def_vis_dis}

Before stating the definition of viscosity solutions for the
discrete problem $(\ref{FK})$-$(\ref{Init_cond_FK})$, we start by
defining some terminology. Let $\eps, \delta>0$ and
$0<T\leq+\infty$. Given a point $P_0=(x_{0},t_{0})\in
\overline{\O}^{\eps}\times [0,T)$, the set ${\mathcal
N}^\varepsilon_\delta(t_0)\subset\overline{\O}^{\eps}\times [0,T)$
is defined by:
\begin{equation}\label{neighborhood}
\mathcal{N}^{\eps}_{\delta}(P_0)=
\{P=(y,t)\in\overline{\O}^{\eps}\times [0,T); |y-x_0| \le
\varepsilon  \mbox{ and } |t-t_{0}|<\delta\}.
\end{equation}
The spaces $USC(\overline{\O}^{\eps}\times [0,T))$ and
$LSC(\overline{\O}^{\eps}\times [0,T))$ are defined respectively by:
\begin{equation}\label{USCt}
USC(\overline{\O}^{\eps}\times [0,T))=\{u\mbox{ defined on
  $\overline{\O}^{\eps}\times [0,T)$};\, \forall x\in
\overline{\O}^{\eps},\, u(x,\cdotp)\in USC([0,T))\}
\end{equation}
and
\begin{equation}\label{LSCt}
LSC(\overline{\O}^{\eps}\times [0,T))=\{u\mbox{ defined on
  $\overline{\O}^{\eps}\times [0,T)$};\, \forall x\in
\overline{\O}^{\eps},\, u(x,\cdotp)\in LSC([0,T))\},
\end{equation}
where $USC([0,T))$ (resp. $LSC([0,T))$) is the set of locally
bounded upper (resp. lower) semi-continuous functions on $[0,T)$. In
a similar manner, we define the space
$C^{k}(\overline{\O}^{\eps}\times [0,T))$, $k\in \N$, by:
\begin{equation}\label{C1t}
C^{k}(\overline{\O}^{\eps}\times [0,T))=\{u\mbox{ defined on
  $\overline{\O}^{\eps}\times [0,T)$};\, \forall x\in
\overline{\O}^{\eps},\, u(x,\cdotp)\in C^{k}([0,T))\}.
\end{equation}
Next, given $\beta\geq0,$ we present the following definition:

{\color{black}
\begin{definition} $(\textbf{Viscosity sub/super
solutions})$\label{definitionsubsuperdis}

\noindent ${\bf1.}$ {\bf Viscosity sub-solutions.} A function $u\in
USC(\overline{\O}^{\eps}\times [0,T))$ is a viscosity
sub-solution of problem $(\ref{FK})$-$(\ref{Init_cond_FK})$ provided that:\\
\noindent ${\bf(i)}$ $u(x,0)\leq u_{0}(x)$ for $\;\left\{
\begin{aligned}
& x\in \overline{\O}^{\eps}\quad&\mbox{if}& \quad \b> 0\\
& x\in \partial \O^{\eps} \quad &\mbox{if}& \quad \beta=0,
\end{aligned}
\right.$\\
\noindent ${\bf(ii)}$ for any $\varphi \in C^{1}(\overline{\O}^{\eps}\times
[0,T))$, if $u- \varphi$ has a zero local maximum at a point
$P_0=(x_{0},t_{0})\in \overline{\O}^{\eps}\times
[0,T)$ $(\mbox{i.e. $\exists\,\delta>0$ such that $\forall P\in
\mathcal{N}^{\eps}_{\delta}(P_0)$, we have $(u-
\varphi)(P) \leq (u - \varphi)(P_0)=0$})$ then
\begin{equation}\label{dis_eq1}
\left\{\begin{array}{ll}
\beta \varphi_{t}(P_0) \leq \Delta^\varepsilon [\varphi](P_0)\,
&
\left|\begin{array}{ll}
\mbox{ for } P_0\in {\O}^{\eps}\times (0,T) & \quad \mbox{if}\quad \beta >0\\
\mbox{ for } P_0\in {\O}^{\eps}\times [0,T) & \quad \mbox{if}\quad \beta =0
\end{array}\right.\\
\varphi_{t}(P_0)\leq F(\varphi(P_0)) + D^{\eps}[\varphi](P_0)\,
&\mbox{ for }\, P_0\in \partial \O^{\eps} \times (0,T).
\end{array}\right.
\end{equation}
\noindent ${\bf2.}$ {\bf Viscosity super-solutions.} A function $u\in
LSC(\overline{\O}^{\eps}\times [0,T))$ is a viscosity
super-solution of problem $(\ref{FK})$-$(\ref{Init_cond_FK})$ provided that:\\
\noindent ${\bf(i)}$ $u(x,0)\geq u_{0}(x)$ for $\;\left\{
\begin{aligned}
& x\in \overline{\O}^{\eps}\quad&\mbox{if}& \quad \b> 0\\
& x\in \partial \O^{\eps} \quad &\mbox{if}& \quad \beta=0,
\end{aligned}
\right.$\\
\noindent ${\bf(ii)}$ for any $\varphi \in C^{1}(\overline{\O}^{\eps}\times
[0,T))$, if $u - \varphi$ has a zero local minimum at a point
$P_0=(x_{0},t_{0})\in \overline{\O}^{\eps}\times
[0,T)$ $(\mbox{i.e. $\exists\,\delta>0$ such that $\forall P\in
\mathcal{N}^{\eps}_{\delta}(P_0)$, we have $(u -
\varphi)(P) \geq (u^{\eps} - \varphi)(P_0)=0$})$ then
\begin{equation}\label{dis_eq2}
\left\{\begin{array}{ll}
\beta \varphi_{t}(P_0) \geq \Delta^\varepsilon [\varphi](P_0)\,
&
\left|\begin{array}{ll}
\mbox{ for } P_0\in {\O}^{\eps}\times (0,T) & \quad \mbox{if}\quad \beta >0\\
\mbox{ for } P_0\in {\O}^{\eps}\times [0,T) & \quad \mbox{if}\quad \beta =0
\end{array}\right.\\
\varphi_{t}(P_0)\geq F(\varphi(P_0)) + D^{\eps}[\varphi](P_0)\,
&\mbox{ for }\, P_0\in \partial \O^{\eps} \times (0,T).
\end{array}\right.
\end{equation}
\noindent ${\bf3.}$ {\bf Viscosity solutions.} A function $u\in
C^{0}(\overline{\O}^{\eps}\times [0,T))$ is a viscosity solution of
problem $(\ref{FK})$-$(\ref{Init_cond_FK})$ if it is a viscosity
sub- and super-solution of $(\ref{FK})$-$(\ref{Init_cond_FK}).$
\end{definition}

\begin{definition}\label{defi2}
We say that $u$ is a viscosity sub-solution $(\mbox{resp.
super-solution})$ of $(\ref{FK})$ if $u$ only satisfies ${\bf (ii)}$
in the Definition $\ref{definitionsubsuperdis}.$
\end{definition}

\begin{rem}
Note that a function $u$ is a viscosity subsolution of $(\ref{FK})$ if it satisfies the viscosity inequalities on
$\overline{\O}^{\eps}\times (0,T)$ if $\beta >0$ and on
$\left(\overline{\O}^{\eps}\times (0,T)\right) \cup \left(\Omega^\varepsilon\times \left\{0\right\}\right)$ if $\beta=0$.
\end{rem}
}

{\color{black}
\begin{theorem}$(\textbf{Comparison principle in the discrete case})$\label{L2}\\
Let $u\in USC(\overline{\O}^\eps\times[0,T))\;(\hbox{resp.}\; v\in
LSC(\overline{\O}^\eps\times[0,T)))$ be a viscosity sub-solution $(\hbox{resp.
super-solution})$ for the problem $(\ref{FK})$-$(\ref{Init_cond_FK})$
where $0<T<\infty,$ such that
$\|u\|_{L^\infty(\overline{\Omega}^\varepsilon\times[0,T))},\|v\|_{L^\infty(\overline{\Omega}^\varepsilon\times[0,T))}<+\infty$.
Then
$$u(x,t)\leq v(x,t) \quad \mbox{for all}\quad
(x,t)\in\overline{\Omega}^\varepsilon\times[0,T).$$
\end{theorem}
Theorem $\ref{L2}$ will be proved in Section~\ref{sec7}. }

\subsection{Viscosity solutions for the continuous problem}\label{def_vis_con}
Similarly as in subsection $\ref{def_vis_dis}$ we denote by
$USC(\overline{\Omega}\times[0,T))$ (resp.
$LSC(\overline{\Omega}\times[0,T))$) the set of locally bounded
upper (resp. lower) semi-continuous functions on
$\overline{\Omega}\times[0,T)$. We also denote by
$C^k(\overline{\Omega}\times [0,T))$ the classical set of
continuously $k$-differentiable functions on
$\overline{\Omega}\times [0,T)$.

{\color{black}
\begin{definition}$(\textbf{Viscosity {\color{black}sub/super} solutions})$\label{definitionsubsupercont}\\
\noindent ${\bf1.}$ {\bf Viscosity sub-solutions.} A function $u\in
USC(\overline{\O}\times [0,T))$ is a viscosity
sub-solution of problem $(\ref{PN})$-$(\ref{Init_cond_PN})$ provided that:\\
\noindent ${\bf(i)}$ $u\leq u_{0}$ on $\left\{\begin{aligned}
&\overline{\O}\times \{0\}\quad&\mbox{if}&\quad\beta>0\\
&\partial\O\times \{0\}\quad&\mbox{if}&\quad\beta=0.
\end{aligned}
\right.$\\
\noindent ${\bf(ii)}$ for any $\varphi \in C^{2}(\overline{\O}\times [0,T))$,
if $u - \varphi$ has a zero local maximum at a point $P_0=(x_{0},t_{0})\in
\overline{\O}\times [0,T)$ then
\begin{equation}\label{viscosity_sub}
\left\{\begin{array}{ll}
\beta \varphi_{t}(P_0)\leq \Delta \varphi(P_0)
&
\left|\begin{array}{ll}
\mbox{ for } P_0\in {\O}\times (0,T) & \quad \mbox{if}\quad \beta >0\\
\mbox{ for } P_0\in {\O}\times [0,T) & \quad \mbox{if}\quad \beta =0
\end{array}\right.\\
\displaystyle \min \left\{\beta \varphi_{t}(P_{0})-\Delta \varphi(P_{0}) \;,\;
\varphi_{t}(P_0) - F(\varphi(P_0)) - \frac{\partial
  \varphi}{\partial x_{n}}(P_0)\right\}\leq 0
&\mbox{ for }\, P_0\in \partial \O \times (0,T).
\end{array}\right.
\end{equation}

\noindent ${\bf2.}$ {\bf Viscosity super-solutions.}  A function
$u\in LSC(\overline{\O}\times [0,T))$ is a viscosity
super-solution of problem $(\ref{PN})$-$(\ref{Init_cond_PN})$ provided that:\\
\noindent ${\bf(i)}$ $u\geq u_{0}$ on $\left\{\begin{aligned}
&\overline{\O}\times \{0\}\quad&\mbox{if}&\quad\beta>0\\
&\partial\O\times \{0\}\quad&\mbox{if}&\quad\beta=0.
\end{aligned}
\right.$\\
\noindent ${\bf(ii)}$ for any $\varphi \in C^{2}(\overline{\O}\times [0,T))$,
if $u - \varphi$ has a zero local minimum at a point $P_0=(x_{0},t_{0})\in
\overline{\O}\times [0,T)$ then
\begin{equation}\label{viscosity_sup}
\left\{\begin{array}{ll}
\beta \varphi_{t}(P_0)\geq \Delta \varphi(P_0)
&
\left|\begin{array}{ll}
\mbox{ for } P_0\in {\O}\times (0,T) & \quad \mbox{if}\quad \beta >0\\
\mbox{ for } P_0\in {\O}\times [0,T) & \quad \mbox{if}\quad \beta =0
\end{array}\right.\\
\displaystyle \max \left\{\beta \varphi_{t}(P_{0})-\Delta \varphi(P_{0}) \;,\;
\varphi_{t}(P_0) - F(\varphi(P_0)) - \frac{\partial
  \varphi}{\partial x_{n}}(P_0)\right\}\geq 0
&\mbox{ for }\, P_0\in \partial \O \times (0,T).
\end{array}\right.
\end{equation}

\noindent ${\bf3.}$ {\bf Viscosity solutions.} A function $u\in
C^0(\overline{\Omega}\times [0,T))$ is a viscosity solution of
problem $(\ref{PN})$-$(\ref{Init_cond_PN})$ if it is a viscosity
sub- and super-solution of $(\ref{PN})$-$(\ref{Init_cond_PN}).$
\end{definition}

\begin{theorem}$(\textbf{Comparison principle in the continuous case})$\label{Comp_pri_con}\\
Let $u\in USC(\overline{\O}\times [0,T))$ $(\mbox{resp. $v\in
LSC(\overline{\O}\times [0,T))$})$ be a viscosity sub-solution
$(\mbox{resp. super-solution})$ of problem
$(\ref{PN})$-$(\ref{Init_cond_PN})$ where $0<T<\infty,$ such that
$\|u\|_{L^\infty(\overline{\Omega}\times[0,T))},\|v\|_{L^\infty(\overline{\Omega}\times[0,T))}<+\infty$.
Then
$$u\leq v \quad \mbox{on}\quad \overline{\Omega}\times[0,T)$$
\end{theorem}
This theorem will be proved in section $\ref{sec8}.$
}

\begin{rem}
Remark that in the case $\beta=0$, the problem can be reformulated
(at least for smooth solutions) as a nonlocal evolution equation
written on the boundary $\partial\Omega.$ See for instance the work
\cite{CS} on the relation bewteen fractional Laplacian and harmonic
extensions. Note that there is also a viscosity theory for nonlocal
operators (see \cite{BI}).
\end{rem}

\section{Construction of barriers}\label{sec6}

This section is devoted to the construction of barriers for the solution
$u^\eps$ of $(\ref{FK})$-$(\ref{Init_cond_FK})$ for all
$\beta\geq0$.

\subsection{Discrete harmonic extension}
Here we construct the discrete harmonic extension of $u_0$ on
$\overline{\O}^\eps$, i.e. we prove the existence of a solution of
the following problem
 \begin{equation}\label{harmonic_discrete}
\left\{
\begin{aligned}
& -\Delta^1[u_0^D]=0,
   &\quad\mbox{in}&\;\Omega^{1}\\
&u_0^D=u_0,
&\quad\mbox{on}&\;\partial{\Omega^{1}},
\end{aligned}
\right.
\end{equation}
where $\Delta^1=\Delta^\eps$ for $\eps=1$ is defined in
$(\ref{deltaeps})$. We note here that, for the sake of simplicity,
we have taken $\eps=1$, while otherwise it can be treated in the
same way (or simply deduced by rescaling). We will say that a function $u$ is a sub-solution (resp.
super-solution) of  $(\ref{harmonic_discrete})$ if the symbol
\textquotedblleft$=$\textquotedblright in
$(\ref{harmonic_discrete})$ is replaced by the symbol
\textquotedblleft$\leq$\textquotedblright (resp.
\textquotedblleft$\geq$\textquotedblright).

{\color{black}
\begin{lem}
If $u_0\in L^\infty(\partial\O^1)$, then there exists a unique
(viscosity) solution of $(\ref{harmonic_discrete})$ on
$\overline{\O}^1$.
\end{lem}}
\noindent\proof $\;$\\
\noindent{\bf Step 1: Existence.}\\
Define
$$\overline{u}=\sup_{\partial\O^1}u_0\quad\mbox{and}\quad \underline{u}=\inf_{\partial\O^1}u_0,$$
then ${\color{black}\underline{u}}$ (resp.
${\color{black}\overline{u}}$) is a sub- (resp. super-) solution of
(\ref{harmonic_discrete}). Consider now the set
$$\mathcal{S}=\left\{u:\overline{\O}^1\rightarrow\R, \;\mbox{subsolution
of}\;(\ref{harmonic_discrete})\;\mbox{s.t.}\;\underline{u}\leq
u\leq \overline{u}\right\},$$ and in order to use Perron's method, we define
$$u_0^D:=\sup_{u\in\mathcal{S}}u.$$
Using the fact that the maximum of two sub-solutions is a
sub-solution, it is possible to show that $u_0^D$ is a sub-solution
of (\ref{harmonic_discrete}). Hence, it remains to prove that
$u_0^D$ is a super-solution for (\ref{harmonic_discrete}). Suppose
that there exists ${\color{black}x_0}\in\overline{\O}^1$ such that
$$
\left\{
\begin{aligned}
& -\Delta^1[u_0^D]({\color{black}x_0})<0,
   &\quad\mbox{if}\;{\color{black}x_0}\in&\,\Omega^{1}\\
   \mbox{or}\\
&u_0^D({\color{black}x_0})<u_0({\color{black}x_0}),
&\quad\mbox{if}\;{\color{black}x_0}\in&\,\partial{\Omega^{1}}.
\end{aligned}
\right.
$$
Then we construct a sub-solution $w\in\mathcal{S}$ such that
$w>u_0^D$ at ${\color{black}x_0}$. Two cases are considered.\\

\noindent {\bf Case ${\bf -\Delta^1[u_0^D]({\color{black}x_0})<0}$
for ${\bf {\color{black}x_0}\in\O^1}$.} Let $w$ be defined by
$$w({\color{black}x})=\left\{\begin{array}{ll}
\displaystyle{u_0^D({\color{black}x})}&\displaystyle{\quad\mbox{if}\quad {\color{black}x}\neq {\color{black}x_0},}\\
\displaystyle{u_0^D({\color{black}x_0})+\frac{1}{2n}\Delta^1[u_0^D]({\color{black}x_0})}&\displaystyle{\quad\mbox{if}\quad
{\color{black}x}= {\color{black}x_0}.}
\end{array}
\right.
$$
We check that $\underline{u}\leq w\leq\overline{u}$, and we compute
$$\Delta^1[w]({\color{black}x})=\left\{\begin{aligned}
&\geq\Delta^1[u_0^D]({\color{black}x})\geq0&\quad\mbox{if}& \quad {\color{black}x}\neq {\color{black}x_0}\\
&0&\quad\mbox{if}&\quad {\color{black}x}= {\color{black}x_0},
\end{aligned}
\right.$$ where we have used the fact that $w\geq u_0^D$ to deduce
the inequality. Therefore, we have $w\in\mathcal{S}$ and $w({\color{black}x_0})>u_0^D({\color{black}x_0})$ which is a contradiction.\\

\noindent {\bf Case ${\bf
u_0^D({\color{black}x_0})<u_0({\color{black}x_0})}$ for ${\bf
{\color{black}x_0}\in\partial\O^1}$.} Let $w$ be defined by
$$w({\color{black}x})=\left\{\begin{array}{ll}
\displaystyle{u_0^D({\color{black}x})}&\displaystyle{\mbox{if}\quad {\color{black}x}\neq {\color{black}x_0},}\\
\displaystyle{u_0({\color{black}x_0})}&\displaystyle{\mbox{if}\quad
{\color{black}x}= {\color{black}x_0}.}
\end{array}
\right.
$$
By construction, we have $\underline{u}\leq w\leq\overline{u}.$
Moreover, as $w\geq u_0^D$ it follows that
$\Delta^1[w]\geq\Delta^1[u_0^D]\geq0$ on $\O^1$. Therefore
$w\in\mathcal{S}$ and
$w({\color{black}x_0})=u_0({\color{black}x_0})>u_0^D({\color{black}x_0});$
this implies a
contradiction.\\

{\color{black}
\noindent{\bf Step 2: Uniqueness.}\\
We simply adapt case 1 i) of the proof of the comparison principle (Theorem \ref{L2}), given in Section \ref{sec7}.
$\hfill\square$
}

\subsection{Continuous harmonic extension}

For the case $\beta=0$,
we need to consider  the harmonic "extension"
$u_0^c:\overline{\Omega}\rightarrow\mathbb{R}$ of the initial data
$u_0$. This function is the solution of
\begin{equation}\label{harmonicextension1}
\left\{\begin{array}{ll}
\displaystyle{-\Delta u_0^c=0}&\displaystyle{\;\mbox{on}\quad\Omega}\\
\displaystyle{ u_0^c=u_0}&\displaystyle{\;\mbox{on}\quad\partial\Omega}.\\
\end{array}
\right.
\end{equation}

This section is devoted to recall the existence of the continuous harmonic extension and to show some of its properties.
Let $z=(z',z_n)\in\O$ where $z'$ is
identified to an element of $\partial\O$ and $z_n>0,$ and let $H$ be
a new function defined by
\begin{equation}\label{poissonkernel}
H(z',z_n):=\frac{2z_n}{\omega_n(z_n^2+{z'}^2)^{n/2}},
\end{equation}
where $\omega_n$ is the measure of the $(n-1)$-dimensional sphere in
${\R}^n.$ Now, we define $u_0^c$ by
\begin{equation}\label{fundamental_solution}
u_0^c(x):=\int_{\partial\O}H(x'-z',x_n)u_0(z')\,dz'\quad\mbox{for
all}\;x=(x',x_n)\in\mathbb{R}^{n-1}\times(0,+\infty).
\end{equation}
Note that $H(x'-z',x_n)$ is the Poisson kernel. We have the
following three lemmas:
\begin{lemma}$(\textbf{Existence of a continuous harmonic
extension})$\label{continuous_harmonic_extension}\\
If $u_0(z')$ is bounded and continuous for $z'\in\partial\O,$
then the function $u_0^c(x)$ defined by
$(\ref{fundamental_solution})$ belongs to $C^\infty(\O)$ and is
harmonic in $\O$ and extends continuously to $\overline{\O}$ such
that $u_0^c=u_0$ on $\partial\O.$
\end{lemma}
\proof $\;$ See \cite[Section 7.3, p. 129]{ABR}.$\hfill\square$
\begin{lemma}$(\textbf{Estimate on the harmonic extension})$\\
If $u_0\in W^{3,\infty}(\partial\O)$, then we have
\begin{equation}\label{properties1}
|Du_0^c(x)|\leq\frac{C}{{\color{black}1+x_n}},\qquad|D^2u_0^c(x)|\leq\frac{C}{1+x_n^2}\quad\mbox{for
all}\;x=(x',x_n)\in\mathbb{R}^{n-1}\times(0,+\infty).
\end{equation}
 \end{lemma}
\noindent {\bf Proof.} Since $u_0$ is bounded, then by adding an
appropriate constant to $u_0$ we can always assume, without loss of
generality, that $u_0\geq0$. From (\ref{fundamental_solution}) and
from the following property (see \cite[ p.126]{ABR})
$$\int_{\partial\O}H(x'-z',x_n)\,dz'=1\quad\mbox{for
all}\;x=(x',x_n)\in\overline{\O},$$
we have:
$$
|u_0^c|\leq
\int_{\partial\O}H(x'-z',x_n)|u_0(z')|\,dz'\leq
\|u_0\|_{L^\infty}
\int_{\partial\O}H(x'-z',x_n)\,dz'=\|u_0\|_{L^\infty},
$$
 and then
 \begin{equation}\label{bornitude}
\|u_0^c\|_{L^\infty}\leq\|u_0\|_{L^\infty}.
\end{equation}
Next, for $x=(x',x_n)\in\mathbb{R}^{n-1}\times(0,+\infty),$ using
$(\ref{poissonkernel})$ and $(\ref{fundamental_solution})$, we have,
on the one hand {\color{black} with $\mu= \frac{2}{\omega_n}$}:
\begin{equation}\label{diff1}
\frac{1}{\mu}\frac{\partial u_0^c}{\partial x_n}=\int_{\partial\O}\frac{u_0(z')}{(x_n^2+(x'-z')^2)^{n/2}}\,dz'-\int_{\partial\O}\frac{nx_n^2u_0(z')}{(x_n^2+(x'-z')^2)^{n/2+1}}\,dz',
\end{equation}
then we obtain
$$\left|\frac{\partial u_0^c}{\partial x_n}\right|\leq\frac{\|{\color{black}u_0}\|_{L^\infty}}{x_n}+n\frac{\|{\color{black}u_0}\|_{L^\infty}}{x_n}\leq\frac{C}{x_n}$$
where we have used $(\ref{bornitude}).$
On the other hand, we have:
\begin{equation}\label{diff2}
\frac{1}{\mu}\frac{\partial u_0^c}{\partial x'}=-n\int_{\partial\O}\frac{x_nu_0(z')(x'-z')}{(x_n^2+(x'-z')^2)^{n/2+1}}\,dz',
\end{equation}
then, using Young's inequality:
\begin{equation}\label{young1}
x_n |x'-z'|\leq \frac{1}{2}(x_n^2+(x'-z')^2),
\end{equation}
we conclude that
$$\left|\frac{\partial u_0^c}{\partial x'}\right|\leq
{\color{black}\frac{n\mu}{2}}
\int_{\partial\O}\frac{u_0(z')}{(x_n^2+(x'-z')^2)^{n/2}}\,dz'{\color{black}\leq\frac{n}{2x_n}}\|{\color{black}u_0}\|_{L^\infty}\leq\frac{C}{x_n}.$$
Therefore
$$|Du_0^c(x)|\leq\frac{C}{x_n}.$$
In order to prove the second inequality in $(\ref{properties1}),$ we differentiate $(\ref{diff1})$ and $(\ref{diff2}),$ we get
$$\frac{1}{\mu}\frac{\partial^2 u_0^c}{\partial x'\partial
x'}=\int_{\partial\O}\frac{u_0(z')[-nx_n{\color{black}I'}]}{(x_n^2+(x'-z')^2)^{n/2+1}}\,dz'+
\int_{\partial\O}\frac{x_n
u_0(z')(n)(n+2)(x'-z')\otimes(x'-z')}{(x_n^2+(x'-z')^2)^{n/2+2}}\,dz',$$
 {\color{black} where $I'$ is the identity matrix of $\R^{n-1}$,} then, as $1/(x_n^2+(x'-z')^2)\leq 1/x_n^2,$ we obtain
$$\frac{1}{\mu}\left|\frac{\partial^2 u_0^c}{\partial x'\partial x'}\right|\leq
n\int_{\partial\O}\frac{u_0(z')x_n}{(x_n^2+(x'-z')^2)^{n/2}x_n^2}\,dz'+
n(n+2)\int_{\partial\O}\frac{u_0(z')x_n(x'-z')^2}{(x_n^2+(x'-z')^2)^{n/2}(x_n^2+(x'-z')^2)^2}\,dz'.$$
Thus, using the following inequality
$$\frac{(x'-z')^2}{(x_n^2+(x'-z')^2)^2}\leq\frac{1}{x_n^2},$$
we deduce that
\begin{equation}\label{estimateZZ}
\left|\frac{\partial^2 u_0^c}{\partial x'\partial x'}\right|\leq\frac{C}{x_n^2}.
\end{equation}
Moreover, we have
\begin{eqnarray*}
{\color{black}\frac{1}{\mu}}\frac{\partial^2 u_0^c}{\partial x'\partial x_n}&=&\int_{\partial\O}\frac{u_0(z')(-n)(x'-z')}{(x_n^2+(x'-z')^2)^{n/2+1}}\,dz'+\int_{\partial\O}\frac{u_0(z')x^2_n(n)(n+2)(x'-z')}{(x_n^2+(x'-z')^2)^{n/2+2}}\,dz',\\
&=&
\frac{1}{x_n^2}\int_{\partial\O}\frac{u_0(z')x_n(-n)x_n(x'-z')}{(x_n^2+(x'-z')^2)^{n/2}(x_n^2+(x'-z')^2)}\,dz'+\int_{\partial\O}\frac{u_0(z')x^2_n(n)(n+2)(x'-z')}{(x_n^2+(x'-z')^2)^{n/2+2}}\,dz'.
\end{eqnarray*}
Then, using $(\ref{bornitude})$ and $(\ref{young1}),$ we infer that
$$\left|\frac{\partial^2 u_0^c}{\partial x'\partial x_n}\right|\leq\frac{C}{x_n^2}.$$
{\color{black}Next} since,
$${\color{black}\frac{1}{\mu}}\frac{\partial^2 u_0^c}{\partial x_n\partial x_n}=-3n\int_{\partial\O}\frac{u_0(z')x_n}{(x_n^2+(x'-z')^2)^{n/2+1}}\,dz'+n(n+2)\int_{\partial\O}\frac{u_0(z')x_n^3}{(x_n^2+(x'-z')^2)^{n/2+2}}\,dz',$$
we use similar arguments as above in order to obtain
$$\left|\frac{\partial^2 u_0^c}{\partial x_n\partial x_n}\right|\leq\frac{C}{x_n^2}.$$
{\color{black}
Notice that this last inequality also follows from (\ref{estimateZZ}) joint to the harmonicity of $u_0^c$.\\
}
Finally, we use Schauder's estimate near the boundary $\partial\O:$
$$|u_0^c|_{C^{2,\alpha}(\overline{B_1(x',0)\cap\O})}\leq C\left(|\Delta u_0^c|_{C^{\alpha}(B_2(x',0)\cap\O)}
+|u_0|_{C^{2,\alpha}(B_2(x',0)\cap\partial\O)}\right)=C|u_0|_{C^{2,\alpha}(B_2(x',0)\cap\partial\O)},$$
for all $x'\in\mathbb{R}^{n-1},$ where $C^{\alpha}$ and
$C^{2,\alpha}$, $\alpha\in(0,1)$, are the H\"older spaces of order
$\alpha$ and $2+\alpha,$ respectively. The ball $B_r(x',0)$,
$r=1,2$, stands for the ball of center $(x',0)$ and radius $r$. We
finally conclude that:
$$\sup_{x'\in\mathbb{R}^{n-1}}|u_0^c|_{C^{2,\alpha}(\overline{B_1(x',0)\cap\O})}
\leq C \sup_{x'\in\mathbb{R}^{n-1}}|u_0|_{C^{2,\alpha}(B_2(x',0)\cap\partial\O)},$$
and the result follows.$\hfill\square$
\begin{lemma}$(\textbf{$\eps$-uniform bound})$\label{estimation_harmonic_discrete}\\
Under the assumption $(\ref{Reg1})$ we have
\begin{equation}\label{properties2}
|D^\eps[u_0^c](x)|\leq \frac{C_1}{{\color{black}1+x_n}},
\qquad|\Delta^\eps [u_0^c](x)|\leq\frac{C_1}{1+x_n^2}\quad\mbox{for
all}\;\;x=(x',x_n)\in\O^\eps,
\end{equation}
where $C_1>0$ is a positive constant independent on $\eps.$
 \end{lemma}
 \proof $\;$ {\color{black}Let us set  $\rho^-_n=0$ and $\rho_i^\alpha=1$ otherwise. Then
using Taylor's expansion, we have  for $x\in\O^\eps$:
 $$
\eps D^\eps[u_0^c](x):=\sum_{\pm}\sum_{i=1}^{n} \rho^{\pm}_i\left(u_0^c(x\pm\eps e_i)-u_0^c(x)\right)\\
 =\sum_{\pm}\sum_{i=1}^{n} \rho^{\pm}_i\int_0^1\,dtDu_0^c(x\pm\eps te_i)(\pm\eps e_i),
$$
 where $e_i=(0,\dots,0,1,0,\dots,0)$ is the unit vector in $\R^n$ with respect to the $i$ component. Then, using $(\ref{properties1})$, we get
$$
   |D^\eps[u_0^c](x)|\leq (2n-1)\sup_{\begin{aligned}
   1\leq &i\leq n\\
   &\pm
   \end{aligned}}\int_0^1\,dt \rho^{\pm}_i |Du_0^c(x\pm\eps te_i) | \\
   \leq\frac{C_1}{1+x_n}.
$$
}
 Then, in order to terminate the proof, it is sufficient to write $\Delta^\eps[u_0^D]$ using Taylor's expansion as
\begin{equation}\label{eq::comp16}
\begin{array}{lll}
\eps^2 \Delta^\eps[u_0^c](x)&:=& \displaystyle \sum_{\pm}\sum_{i=1}^{n}\left(u_0^c(x\pm\eps e_i)-u_0^c(x)\right)\\
\\
 &=&\displaystyle \sum_{\pm}\sum_{i=1}^{n}\left(\pm\eps e_i\nabla u_0^c(x)+\int_0^1\,dt\int_0^t\,ds D^2u_0^c(x\pm\eps se_i)(\pm\eps e_i)(\pm\eps
 e_i)\right)\\
\\
 &=&\displaystyle \eps^2\sum_\pm\sum_{i=1}^{n}\int_0^1\,dt\int_0^t\,ds D^2u_0^c(x\pm\eps
 se_i)\cdot (e_i,e_i).
\end{array}
\end{equation}
Thus, thanks to (\ref{properties1}), we finally obtain:
$$
   |\Delta^\eps[u_0^c](x)|\leq2n\sup_{\begin{aligned}
   1\leq &i\leq n\\
   &\pm
   \end{aligned}}\int_0^1\,dt\int_0^t\,ds |D^2u_0^c(x\pm\eps se_i) | \\
   \leq\frac{C_1}{1+x_n^2},
   $$
 and the proof is done. $\hfill\square$

\subsection{Uniform barriers for $\beta\ge 0$}\label{barriers2}

{\color{black}
  In this subsection, we show uniform barriers for the solution $u^\eps$ of $(\ref{FK})$-$(\ref{Init_cond_FK})$ for all $\beta\geq0$ in the special case where $u_0=u_0^c$.
\begin{proposition}$(\textbf{Uniform barriers in $\eps$ for all $\beta \geq 0$ for $u_0=u_0^c$})$\label{barriers}\\
Under assumption $(\ref{Reg1}),$ there exists a constant $C>0$
independent on $\eps>0,$ $\beta\ge 0$ and $T$ such that if
\begin{equation}\left\{\begin{array}{l}
\overline{u}^+(x,t):=u_0^c(x)+C(\sqrt{1+x_n}-1+t)\\
\overline{u}^-(x,t):=u_0^c(x)-C(\sqrt{1+x_n}-1-t)\\
\end{array} \right| \quad \mbox{for}\quad (x,t)\in\overline{\O}^{\eps}\times[0,T),
\end{equation}
then $\overline{u}^+$ (resp. $\overline{u}^-$) is a supersolution
(resp. subsolution) of $(\ref{FK})$-$(\ref{Init_cond_FK})$. However,
if $u^\varepsilon$ is a bounded viscosity solution of
$(\ref{FK})$-$(\ref{Init_cond_FK})$, then we have
\begin{equation}
\overline{u}^-\leq u^\eps\leq
\overline{u}^+\quad\mbox{on}\;\;\overline{\O}^\eps\times[0,T),
\end{equation}
and moreover:
$$|u^\varepsilon (x,t)|\le {\color{black}\|u_0\|}_{L^\infty} + t {\color{black}\|F\|}_{L^\infty}
\quad \mbox{for}\quad (x,t)\in\overline{\O}^{\eps}\times[0,T).$$
\end{proposition}
\noindent {\bf Proof.}\\
\noindent {\bf Step 1: sub/supersolution property}\\
We first check that $\overline{u}^+$ is a super-solution of $(\ref{FK})$-$(\ref{Init_cond_FK})$.
Indeed,\\
\noindent $\bullet\;$ If $(x,t)\in\overline{\O}^\eps\times\{0\},$ then $\overline{u}^+(x,t)\ge u_0^c(x)=u_0(x).$\\
\noindent $\bullet\;$ If $(x,t)\in\partial\O^\eps\times(0,T),$ then we take $C>0$ such that
$$C\geq  F(\overline{u}^+(x,t))+D^\eps[u_0^c](x)+CD^\eps [\sqrt{1+x_n}]\quad\Longleftrightarrow\quad \overline{u}^+_t(x,t)\geq
F(\overline{u}^+(x,t))+D^\eps[\overline{u}^+](x,t).$$
Using the fact that $D^\eps [\sqrt{1+x_n}](x_n=0)\le \frac12$, we see that such a constant $C$ exists.\\
\noindent $\bullet\;$ If $(x,t)\in\O^\eps\times[0,T),$ then we can use Lemma \ref{estimation_harmonic_discrete}.
We get
$$\Delta^\eps [u_0^c](x)\leq \frac{C_1}{1+x_n^2}\leq\frac{C}{4(1+x_n)^{3/2}},$$
where the last inequality is true for $C>0$ large enough. Moreover,
by repeating the same computations as in (\ref{eq::comp16}), we
easily get with $f(a)=\sqrt{1+a}$:
\begin{equation}\label{eq::concave}
\Delta^\eps [\sqrt{1+x_n}] ={\color{black}\int_0^1dt\int_0^t\,ds
\sum_{\pm} f''(x_n\pm \varepsilon s)} \le{\color{black} \int_0^1
dt\int_0^t\,ds.2. f''(x_n)} \leq f''(x_n)=-\frac{1}{4(1+x_n)^{3/2}},
\end{equation}
where we have used the concavity of the function $f''$.
Then we have
 $$\Delta^\eps[u_0^c](x)\leq \frac{C}{4(1+x_n)^{3/2}}\leq -C\Delta^\eps[\sqrt{1+x_n}],$$
 hence
 $$\beta \overline{u}_t^+(x,t)=\beta C\geq0\geq\Delta^\eps [\overline{u}^+](x,t).$$
As a consequence, we finally deduce that $\overline{u}^+$ is a
super-solution. In a similar way, by taking in addition the
following assumption on $C:$
$$C\geq  -F(\overline{u}^-(x,t))-D^\eps[u_0^c](x)+CD^\eps [\sqrt{1+x_n}]$$
we can prove that $\overline{u}^-$ is a sub-solution.\\

\noindent {\bf Step 2: bounds on $u^\varepsilon$}\\
Let us call $\bar{v}(x,t):={\color{black}\|u_0\|}_{L^\infty} + t
{\color{black}\|F\|}_{L^\infty}$. It is easy to check that $\bar{v}$
is a super{\color{black}-}solution of
$(\ref{FK})$-$(\ref{Init_cond_FK})$. Then
$$\min(\overline{u}^+,\bar{v})$$
is still a super{\color{black}-}solution and is bounded. Therefore,
if $u^\varepsilon$ is a bounded viscosity solution of
$(\ref{FK})$-$(\ref{Init_cond_FK})$, we can then apply the
comparison principle (Theorem \ref{L2}) to conclude that
$$u^\varepsilon \le \min(\overline{u}^+,\bar{v}).$$
This shows the upper bounds. For the lower bounds, we proceed similarly with $\max(\overline{u}^-,-\bar{v})$.
$\hfill\square$
}

\subsection{Barriers for $\beta> 0$}
{\color{black} We have the following
\begin{proposition}$(\textbf{Barriers for $\beta >0$})$\label{barriers+}\\
Under the assumptions of Theorem \ref{theo1}, for every $\beta >0$,
there exists a constant $C_\beta >0$ such that for all
$\varepsilon$, if $u^\varepsilon$ is a bounded viscosity
 solution of $(\ref{FK})$-$(\ref{Init_cond_FK})$, then we have
$$|u^\varepsilon (x,t)-u_0(x)|\le t C_\beta
\quad \mbox{for}\quad (x,t)\in\overline{\O}^{\eps}\times[0,T).$$
\end{proposition}

\noindent {\bf Proof.}\\
Using in particular the fact that $F$ is bounded and that $u_0\in
W^{2,\infty}(\overline{\Omega})$, we simply check that $u_0(x) +
tC_\beta$ is a super{\color{black}-}solution for a suitable constant
$C_\beta{\color{black}>0}$ and apply the comparison principle. We
proceed similarly for sub{\color{black}-}solutions $u_0(x) -
tC_\beta$.$\hfill{\square}$}

\section{Existence and uniqueness of a solution for the discrete problem}\label{sec5}

The aim of this section is to prove the existence of solutions of
problem $(\ref{FK})$-$(\ref{Init_cond_FK}).$ Cauchy-Lipschitz method
is the main tool used to prove the existence of solutions for
$\beta>0,$ while in the case $\beta=0$ we need barriers to prove the
existence.

\begin{thm}$(\textbf{Existence and uniqueness, $\beta>0$})$\label{H}\\
If $u_0$ and $F$ satisfy $(\ref{Reg1})$ and $\beta>0$, then there
exists a unique bounded solution $u^{\beta,\varepsilon}\in
C^1(\overline{\Omega}^\varepsilon\times[0,T))$ of problem
$(\ref{FK})$-$(\ref{Init_cond_FK}).$
\end{thm}
\proof $\;$ The proof is done using the classical Cauchy-Lipschitz
theorem. Let
${\color{black}B:=L^\infty(\overline{\Omega}^\varepsilon)}$ be the
Banach space with the norm
$${\color{black}\|u\|_B:=\sup_{x\in\overline{\Omega}^\varepsilon}|u(x)|,}\qquad\mbox{for
every}\;\;u\in B,$$ and $\mathcal{F}:B\longrightarrow B$ be the map
defined, for every $u\in B$ and
${\color{black}x\in\overline{\O}^\eps},$ by
$$\mathcal{F}[u]({\color{black}x}):=\left\{
\begin{array}{ll}
\displaystyle {\frac{1}{\beta}\Delta^\eps[u]({\color{black}x}),}&\displaystyle{{\color{black}x\;\in\Omega^\varepsilon},}\\
   \displaystyle {F(u({\color{black}x}))+D^\eps[u]({\color{black}x}),}&\displaystyle{{\color{black}x\;\in\partial\Omega^\varepsilon}}
\end{array}
\right.
$$
{\color{black}where $\Delta^\varepsilon[u](x)$ and
$D^\varepsilon[u](x)$ are defined as in (\ref{deltaeps}), dropping
the variable $t$ on both sides of (1.6).} Then, for every $u,v\in B$
and ${\color{black}x\in\overline{\Omega}^\varepsilon},$ we have two
cases; either ${\color{black}x\in\Omega^\varepsilon}$, and hence we
obtain
\begin{eqnarray*}
  \|\mathcal{F}[u]({\color{black}x})-\mathcal{F}[v]({\color{black}x})\|_B &\leq&\|\frac{1}{\beta}\Delta^\eps[u-v]({\color{black}x}) \|_B\\
   &\leq&\frac{4n}{\beta\varepsilon^2}\|u-v\|_B,
\end{eqnarray*}
or ${\color{black}x\in\partial\Omega^\varepsilon}$, then:
\begin{eqnarray*}
  \|\mathcal{F}[u]({\color{black}x})-\mathcal{F}[v]({\color{black}x})\|_B &\leq&\|D^\eps[u-v]({\color{black}x})+|F(u({\color{black}x}))-F(v({\color{black}x}))|\|_B\\
  &\leq&\left(\frac{2(2n-1)}{\varepsilon}+\|F'\|_{L^\infty(\mathbb{R})}\right)\|u-v\|_B.
\end{eqnarray*}
In all cases we conclude that $\mathcal{F}$ is globally Lipschitz
continuous. By the Cauchy-Lipschitz theorem, we get the existence
and uniqueness of a solution $u^{\beta,\varepsilon}\in
{\color{black}C^1(\overline{\Omega}^\varepsilon\times[0,T))}$
satisfying
$${\color{black}u^{\beta,\varepsilon}_t(\cdot,t)={\mathcal F}[u^{\beta,\varepsilon}](\cdot,t),}$$
{\color{black} such that
$u^{\beta,\varepsilon}(\cdotp,0)=u_0(\cdotp).$} $\hfill\square$

\begin{thm}$(\textbf{Existence and uniqueness, $\beta=0$})$\label{H0}\\
If $u_0$ and $F$ satisfy $(\ref{Reg1})$ and $\beta=0$, then there
exists a unique bounded continuous solution $u^{0,\eps}$ on
$\overline{\Omega}^\varepsilon\times[0,T)$ of the problem
$(\ref{FK})$-$(\ref{Init_cond_FK}).$
\end{thm}
{\color{black} \proof $\;$ We consider the solution
$u^{\beta,\varepsilon}$ given by Theorem \ref{H} for the choice
$u_0=u_0^c$. Let
$$\hat{u}=\limsup_{\beta\rightarrow0}{}^*u^{\beta,\varepsilon}\quad\mbox{and}\quad
\check{u}=\liminf_{\beta\rightarrow0}{}_*u^{\beta,\varepsilon}.$$
Then, using Propositon $\ref{barriers}$, we obtain
\begin{equation}\label{eq::bar}
\overline{u}^-\leq\check{u}\leq\hat{u}\leq \overline{u}^+ \quad
\mbox{and}\quad  |\check{u}|, |\hat{u}|\le
{\color{black}\|u_0\|}_{L^\infty} +
t{\color{black}\|F\|}_{L^\infty}.
\end{equation}
Using standard arguments similar to those in Step 1.1 of the proof
of Theorem \ref{theo1} in Section \ref{sec4}, we can show that
$\check{u}$ (resp. $\hat{u}$) is a super- (resp. sub-) solution of
the problem $(\ref{FK})$. The only difficulty is to recover the
viscosity inequality on $\Omega^\varepsilon \times
\left\{0\right\}$. This last inequality follows from the fact that
$u^{\beta,\varepsilon}$ is a classical solution and then satisfies
the equation also at $t=0$. Finally, using (\ref{eq::bar}), we
deduce that $\check{u}$ (resp. $\hat{u}$) is a super- (resp. sub-)
solution of the problem $(\ref{FK})$-$(\ref{Init_cond_FK}).$ Then
the comparison principle (Theorem~\ref{L2}) implies that
$\hat{u}\leq\check{u}$. So $\check{u}=\hat{u}=:u^{0,\varepsilon}$ is
a continuous bounded solution of
$(\ref{FK})$-$(\ref{Init_cond_FK})$. The uniqueness of this solution
follows again from the comparison principle.
$\hfill\square$\\
}

As a conclusion, there exists a unique (viscosity) solution $u^\varepsilon\in C^0(\overline{\Omega}^\varepsilon\times[0,T))$ of problem
$(\ref{FK})$-$(\ref{Init_cond_FK})$ defined by:
\begin{equation}\label{uepsilon}
u^\varepsilon:=\left\{
\begin{array}{ll}
\displaystyle {u^{\beta,\varepsilon}}&\displaystyle{\mbox{if}\;\;\beta>0,}\\
{}\\
   \displaystyle {u^{0,\varepsilon}}&\displaystyle{\mbox{if}\;\;\beta=0.}
\end{array}
\right.
\end{equation}
\begin{rem}
The existence of a solution could also be proven by Perron's method.
\end{rem}

\section{Proof of Theorem \ref{theo1}}\label{sec4}

This section is devoted to the proof of Theorem $\ref{theo1}.$
Note that
in the sequel, we use the following notation:
$$\O_T:=\O\times(0,T),\quad{\color{black} \partial^l\Omega_T}:=\partial\O\times(0,T),\quad\overline{\O}_T:=\overline{\O}\times[0,T),$$
and
$$\O^\eps_T:=\O^\eps\times(0,T),\quad{\color{black}\partial^l\O^\eps_T}:=\partial\O^\eps\times(0,T),\quad\overline{\O}^\eps_T:=\overline{\O}^\eps\times[0,T),$$
{\color{black}where $\partial^l$ denotes the lateral boundary.}\\

\noindent \textbf{Proof of Theorem~\ref{theo1}.} The proof is
divided into two steps.\\

\noindent{\bf \small{Step 1: }$\;$${\bf\overline{u}}$ and
${\bf\underline{u}}$
are, respectively, sub- and super-solution of ${\bf(\ref{PN})}$-${\bf(\ref{Init_cond_PN}).}$}\\
Let $u^\varepsilon$ be the bounded viscosity solution of the problem
$(\ref{FK})$-$(\ref{Init_cond_FK})$. We define, for
$(x,t)\in{\color{black}\overline{\Omega}_T},$ the functions
$\overline{u}$ and $\underline{u}$ as follows:
\begin{equation}\label{limsupliminf}
\overline{u}(x,t):=\limsup_{{\small\begin{aligned}
y\rightarrow x,s\rightarrow t,\varepsilon\rightarrow0\\
y\in\overline{\O}^\varepsilon,\;s\in[0,T)
\end{aligned}}}
u^\varepsilon(y,s)\qquad\mbox{and}\qquad\underline{u}(x,t):=\liminf_{{\small\begin{aligned}
y\rightarrow x,s\rightarrow t,\varepsilon\rightarrow0\\
y\in\overline{\O}^\varepsilon,\;s\in[0,T)
\end{aligned}}} u^\varepsilon(y,s).
\end{equation}
We start by showing that $\overline{u}$ is a viscosity sub-solution of
$(\ref{PN})$-$(\ref{Init_cond_PN}).$\\

\noindent{\bf \small{Step 1.1: }$\;$Proof of 1 (ii) in
{\color{black}D}efinition
${\bf\ref{definitionsubsupercont}.}$}\\
Let $\varphi\in C^2({\color{black}\overline{\Omega}_T})$ and
$P_0=(x_0,t_0)\in\overline{\Omega}_T$ such that
$\overline{u}-\varphi$ has a {\color{black}zero} local maximum at
$P_0.$ Without loss of generality, we can assume that this maximum
is global and strict. Therefore we have
\begin{equation}\label{P_10}
\forall r>0,\;\exists \delta=\delta(r)>0\;\;\mbox{such
that}\;\varphi-\overline{u}\geq\delta>0\;\;\mbox{in}\;\overline{\Omega}_T\backslash
B_r(P_0),
\end{equation}
where $B_r(P_0)\subset\mathbb{R}^n\times\mathbb{R}$ is the open ball
of radius $r$ and of center $P_0$. Now, let
$w^\varepsilon:=\varphi-u^\varepsilon$ for some $\varepsilon>0.$
Then, using the definition of $\overline{u},$ and inequality
$(\ref{P_10}),$ we infer that
\begin{equation}\label{P_11}
w^\varepsilon\geq\frac{\delta}{2}\qquad\mbox{in}\;\;\overline{\Omega}_T^\varepsilon\backslash
B_r(P_0),
\end{equation}
for all $r>0$ and for $\varepsilon>0$ small enough. Using
\cite[Lemma $4.2$]{Bar}, it is then classical to see that there
exists a sequence
$P^*_\varepsilon=(x^*_\varepsilon,t^*_\varepsilon)\in
\overline{\Omega}_T^\varepsilon\cap B_r(P_0)$ such that, as
$\varepsilon\rightarrow0,$ we have
\begin{equation}\label{local}
P^*_\varepsilon\rightarrow
P_0,\;\;u^\varepsilon(P^*_\varepsilon)\rightarrow
u(P_0)\;\;\hbox{and}\;\;u^\varepsilon-\varphi\;\;\hbox{has a local maximum at $P^*_\varepsilon.$}
\end{equation}
Two cases are then considered.\\

\noindent{\bf\small{Case ${\bf1:}$}}$\;$ $P_0\in\Omega \times
[0,T).$ We argue by contradiction. Assume that there exists a
positive constant $\gamma>0$ such that {\color{black}
\begin{equation}\label{contr1}
\beta\varphi_t(P_0)=\Delta\varphi(P_0)+\gamma
\quad \left|\begin{array}{l}
\mbox{for}\quad P_0\in \Omega\times (0,T) \quad \mbox{if}\quad \beta >0\\
\mbox{for}\quad P_0\in \Omega\times [0,T) \quad \mbox{if}\quad \beta =0\\
\end{array}\right.
\end{equation}
}
Moreover, using $(\ref{local})$ and the fact that $u^\varepsilon$ is
a viscosity sub-solution of problem
$(\ref{FK})$-$(\ref{Init_cond_FK})$, we conclude that
\begin{equation}\label{P_12}
\beta\varphi_t(P^*_\varepsilon)\leq\Delta^\varepsilon[\varphi](P_\varepsilon^*),
\end{equation}
and by Taylor's expansion, this inequality $(\ref{P_12})$ implies
\begin{equation}\label{contr2}
\beta\varphi_t(P^*_\varepsilon)\leq\Delta\varphi(P^*_\varepsilon)+{\color{black}o_\varepsilon(1)}.
\end{equation}
Combining $(\ref{contr1})$ and $(\ref{contr2})$ yields:
$$\gamma\leq-\beta(\varphi_t(P^*_\varepsilon)-\varphi_t(P_0))+\Delta\varphi(P^*_\varepsilon)-\Delta\varphi(P_0)
+o_\varepsilon(1),$$
where the right hand side goes to zero as $\varepsilon\rightarrow0.$
This contradicts the fact that $\gamma>0.$\\

\noindent{\bf \small{Case ${\bf2:}$}}$\;$
$P_0\in{\color{black}\partial^l\Omega_T}.$ We repeat similar
arguments as in Case $1.$ Suppose that
\begin{equation}\label{P_13}
\min\left\{\beta\varphi_t(P_0)-\Delta\varphi(P_0)\;\;,\;\;\varphi_t(P_0)-F(\varphi(P_0))-\frac{\partial \varphi}{\partial
  x_n}(P_0)\right\}=\gamma_1>0.
\end{equation}
Inequality $(\ref{P_13})$ implies
\begin{equation}\label{contr4}
\beta\varphi_t(P_0)-\Delta\varphi(P_0)\geq\gamma_1,
\end{equation}
and
\begin{equation}\label{contr5}
\varphi_t(P_0)-F(\varphi(P_0))-\frac{\partial \varphi}{\partial
  x_n}(P_0)\geq\gamma_1.
\end{equation}
On the one hand, if $P^*_\eps\in\Omega_T^\varepsilon,$ then, using
$(\ref{contr2})$ and $(\ref{contr4}),$ we obtain a contradiction by
using the same reasoning as in Case 1. On the other hand, if
$P^*_\eps\in{\color{black}\partial^l\Omega_T^\varepsilon}$, then
using $(\ref{local})$ and the fact that $u^\varepsilon$ is a
subsolution of $(\ref{FK})$-$(\ref{Init_cond_FK})$ we obtain:
$$\varphi_t(P^*_\varepsilon)\leq F(\varphi(P^*_\eps))+D^\eps[\varphi](P^*_\eps),$$
and then, using Taylor's expansion, we obtain
$$\varphi_t(P^*_\varepsilon)\leq\frac{\partial\varphi}{\partial x_n}(P^*_\varepsilon)+F(\varphi(P^*_\eps))+O(\varepsilon).$$
Finally, subtracting this inequality from $(\ref{contr5})$, we
conclude, after passing to the limit as $\varepsilon\rightarrow0,$
that $\gamma_1\leq0$; contradiction.\\

\noindent{\bf \small{Step ${\bf1.2:}$ }Proof of ${\bf1 (i)}$ in Definition ${\bf\ref{definitionsubsupercont}.}$}\\
\noindent From Propositions \ref{barriers} and \ref{barriers+}, we can pass to the
limit and get from the barriers
 $$|\overline{u}(x,t)-u_0(x)|\leq \left\{\begin{array}{ll}
t C_\beta & \quad \mbox{if}\quad \beta >0\\
C(t+\sqrt{1+x_n}-1) & \quad \mbox{if}\quad \beta =0 \quad \mbox{and}\quad u_0=u_0^c
 \end{array}\right|\quad\mbox{for all}\quad(x,t)\in{\color{black}\overline{\O}_T}.$$
 Therefore for $t=0$ we recover $1 {\bf(i)}$ in Definition $\ref{definitionsubsupercont}$ for all
 $\beta\geq0$.\\

\noindent{\bf \small{Step 2: }Existence and convergence.}\\
\noindent From Step 1 we conclude that $\overline{u}$ is a viscosity
sub-solution of $(\ref{PN})$-$(\ref{Init_cond_PN})$ and by a similar
manner, we can show that $\underline{u}$ is a viscosity
super-solution of $(\ref{PN})$-$(\ref{Init_cond_PN}).$
{\color{black} Moreover $\overline{u}$ and $\underline{u}$ are
bounded, because Proposition \ref{barriers} implies
$$|\overline{u}|,|\underline{u}|\le {\color{black}\|u_0\|}_{L^\infty} + t {\color{black}\|F\|}_{L^\infty}$$
} Then by the comparison principle for problem
$(\ref{PN})$-$(\ref{Init_cond_PN})$ (Theorem $\ref{Comp_pri_con}$),
we have $\overline{u}\leq\underline{u}$. On the other hand, by the
definition of $\overline{u}$ and $\underline{u},$ we have
$\underline{u}\leq\overline{u}.$ As a consequence, we deduce, for
$(x,t)\in{\color{black}\overline{\Omega}_T},$ that:
\begin{equation}\label{equality1}
\overline{u}(x,t)=\underline{u}(x,t)=\lim_{\varepsilon\rightarrow0,y\rightarrow
x,s\rightarrow t}u^\varepsilon (y,s)=:u^0(x,t)
\end{equation}
where $u^0$ is then a continuous viscosity solution of problem
$(\ref{PN})$-$(\ref{Init_cond_PN}).$ Moreover $u^0$ is unique, still
by the comparison principle. Finally,
we also observe that the convergence $u^\varepsilon\rightarrow u^0$
as $\varepsilon\rightarrow0$ is locally uniform in the following
sense:
$${\color{black}\|}u^\varepsilon-u^0{\color{black}\|}_{L^\infty(K\cap\;{\color{black}\overline{\Omega}_T^\varepsilon})}\longrightarrow0\qquad\mbox{as}\qquad
\varepsilon\rightarrow0,$$ for any compact set
$K\subset{\color{black}\overline{\Omega}_T}.$ Indeed, suppose that
there exists $\theta>0$ such that for all $k>0$ there exists
$\eps_k$ satisfying $0<\varepsilon_k<\frac{1}{k}$ with
$$ {\color{black}\|}u^{\varepsilon_k}-u^0 {\color{black}\|}_{L^\infty(K\cap(\overline{\Omega}^{\varepsilon_k}\times[0,T)))}>\theta.$$
Then, there exists a sequence $P_k\in
K\cap(\overline{\Omega}^{\varepsilon_k}\times[0,T))$ such that
\begin{equation}\label{P_24}
|u^\varepsilon(P_k)-u^0(P_k)|>\theta.
\end{equation}
Now, since $K\cap(\overline{\Omega}^{\varepsilon_k}\times[0,T))\subset K$ with
$K$ compact, there exists a subsequence, denoted for simplicity by
$P_k$, such that $P_k\rightarrow P_\infty$ as $k\rightarrow\infty.$ Finally, taking the
$\liminf_{k\rightarrow\infty, P_k\rightarrow P_\infty}$ in the
inequality $(\ref{P_24})$ and using $(\ref{equality1}),$ we obtain
$0>\theta$ which gives a contradiction. $\hfill\square$

\subsection{Other barriers and comments on another possible approach for $\beta=0$}

Let us notice that for $\beta=0$, we have natural barriers for the $\varepsilon$-problem, which are:
$$u_0^{D,\varepsilon}(x) \pm C_\varepsilon t$$
where $u_0^{D,\varepsilon}=u_0^{D}$ is the discrete harmonic extension (we make its dependence explicit on $\varepsilon$)
associated to the operator $\Delta^\varepsilon$, and with the constant $C_\varepsilon$ satisfying
$$C_\varepsilon \ge D^\varepsilon[u_0^{D,\varepsilon}] + {\color{black}\|}F{\color{black}\|}_{L^\infty}.$$
Then another possible approach to show the convergence of $u^\varepsilon$ to $u^0$,
could be to control the solution as $\varepsilon$ goes to zero, showing that:\\
1) the constant $C_\varepsilon$ can be taken independent on $\varepsilon$ (using the regularity of $u_0$ on $\partial\Omega$).\\
2) the discrete harmonic extension $u_0^{D,\varepsilon}$ converges to the continuous harmonic extension $u_0^c$ as $\varepsilon$ tends to zero.\\
Then we could also introduce a (more classical) notion of viscosity
solution in the case $\beta=0$, assuming that at $t=0$, we can
compare the sub/super{\color{black}-}solution to the initial data
(taken to be equal to the harmonic extension $u_0^{D,\varepsilon}$
for the $\varepsilon$-problem and $u_0^c$ for the limit problem).
Nevertheless, this other approach would require some additional work
(to show 1) and 2)), and would not simplify the proofs.

\begin{rem}{\bf (Convergence of the discrete harmonic extension to the continuous harmonic extension)}\\
Notice that point $2)$ is a consequence of our convergence theorem $(\hbox{Theorem $\ref{theo1}$})$ in the case $\beta=0$.
Indeed, in the case $\beta=0$, the bounded solutions satisfy
$$|u^\varepsilon(x,t) - u_0^{D,\varepsilon}(x)|\le C_\varepsilon t$$
and then
$$u^\varepsilon(x,0)=u_0^{D,\varepsilon}(x).$$
On the other hand, the functions
$$u_0^c(x) \pm C_0 t$$
are barriers for the limit problem if
$$C_0\ge  {\color{black}\|}\frac{\partial u_0^c}{\partial n} {\color{black}\|}_{L^\infty}+  {\color{black}\|}F {\color{black}\|}_{L^\infty}.$$
Then bounded solutions of the continuous problem satisfy
$$|u^0(x,t)-u_0^c(x)|\le C_0 t$$
and then
$$u^0(x,0)=u_0^c(x).$$
Finally from the $(\hbox{locally uniform})$ convergence of
$u^\varepsilon$ to $u^0$ in particular for $t=0$, we deduce that
$u_0^{D,\varepsilon}$ converges, locally uniformly, to $u_0^c$.
\end{rem}

\section{Proof of Theorem ${\bf\ref{L2}}$}\label{sec7}

In order to emphasize the main points, we perform the proof in several steps.\\

\noindent{\bf \small{Step 1:} Rescaling.} We want to reduce the
problem $(\ref{FK})$-$(\ref{Init_cond_FK})$
 to the case $\varepsilon=1$ and with a nonlinearity $F$ replaced by a monotone one. To this end, we introduce the new functions
$$\overline{u}(x,t):=e^{-\lambda t}u(\varepsilon x,\varepsilon t),\quad
\overline{v}(x,t):=e^{-\lambda t}v(\varepsilon x,\varepsilon
t),\quad (x,t)\in \overline{\Omega}^1\times [0,\overline{T}),$$
where $\overline{T}=T/\varepsilon$ and $\lambda>0$ is a constant to
be determined later. We see easily that $\overline{u}$
$(\hbox{resp.}\;\overline{v})$ is a sub-solution (resp.
super-solution) of the following problem
\begin{equation}\label{newFK}\left\{
  \begin{array}{ll}
   \,\,\displaystyle{\overline{\beta}\overline{u}_t=\Delta^1[\overline{u}]-\overline{\beta}\lambda\overline{u}}&\displaystyle{\quad\hbox{in}\quad\Omega^1\times(0,\overline{T})},\\
   {}\\
   \displaystyle{\overline{u}_t=\overline{F}(\overline{u},t)+D^1[\overline{u}]}&\displaystyle{\quad\hbox{on}\quad\partial\Omega^1\times(0,\overline{T})},
  \end{array}
\right.
\end{equation}
and
\begin{equation}\label{newInit_cond_FK}
\overline{u}(x,0)=u_{0}(\varepsilon x)\quad \mbox{for} \quad
\left\{
\begin{aligned}
&x\in \overline{\O}^{1}\quad
&\mbox{if}& \quad \overline{\b}> 0\\
& x\in \partial \O^{1} \quad &\mbox{if}& \quad \overline{\b}=0,
\end{aligned}
\right.
\end{equation}
with
$$\overline{\b}=\eps\b\quad\mbox{and} \quad\overline{F}(\overline{u},t)=\eps e^{-\lambda t}F(e^{\lambda t}\overline{u})-\lambda \overline{u}.$$
We argue by contradiction assuming that
$$M:=\sup_{\overline{\Omega}^1\times [0,\overline{T})}(\overline{u}-\overline{v})>0.$$
From the definition of $M,$ there exists a sequence
$P^k=(x^k,t^k)\in\overline{\O}^1\times[0,\overline{T})$ such that
$\overline{u}(P^k)-\overline{v}(P^k)\rightarrow M>0$ as
$k\rightarrow\infty.$ Let us choose an index $k_0$ such that
$$\overline{u}(P^{k_0})-\overline{v}(P^{k_0})\geq \frac{M}{2}.$$
At this stage, if we write $x=(x',x_n)$ with $x'=(x_1,\dots,x_{n-1}),$ we define
$$\Psi(x,t,s):=\frac{|t-s|^2}{2\delta}+\frac{\eta}{\overline{T}-t}+\alpha|x'|^2+\gamma\sqrt{1+x_n},\quad (x,t,s)\in\overline{\O}^1\times[0,\overline{T})^2,$$
where $\delta,\eta,\alpha,\gamma>0$ will be chosen later, and we consider
\begin{equation}\label{P_18}
\overline{M}:=M_{\delta,\eta,\alpha,\gamma}:=\sup_{x\in\overline{\Omega}^1, t,s\in[0,\overline{T})}(\overline{u}(x,t)-\overline{v}(x,s)-\Psi(x,t,s)).
\end{equation}

{\color{black} \noindent{\bf \small{Step 2:} A priori estimates.}
Now, if we choose $\eta,\alpha,\gamma$ such that:
$$\frac{\eta}{\overline{T}-t^{k_0}}+\alpha|(x^{k_0})'|^2+\gamma\sqrt{1+x^{k_0}_n}\leq\frac{M}{4},$$
we conclude that
$$\overline{M}\geq\overline{u}(P^{k_0})-\overline{v}(P^{k_0})-\Psi(x^{k_0},t^{k_0},t^{k_0})\geq\frac{M}{4}>0.$$
Moreover, from the definition of $\overline{M}$, there exists
$(\overline{x},\overline{t},\overline{s})\in\overline{\O}^1\times[0,\overline{T})^2$ such that
\begin{equation}\label{P_28}
\overline{u}(\overline{x},\overline{t})-\overline{v}(\overline{x},\overline{s})-\Psi(\overline{x},\overline{t},\overline{s})=\overline{M}>0
\end{equation}
and
\begin{equation}\label{P_29}
\frac{|\overline{t}-\overline{s}|^2}{2\delta}+\frac{\eta}{\overline{T}-\overline{t}}+\alpha|\overline{x}'|^2+\gamma\sqrt{1+\overline{x}_n}\leq
C
\end{equation}
where
$$C=\|\overline{u}\|_{L^{\infty}(\overline{\O}^1\times[0,\overline{T}))}+\|\overline{v}\|_{L^{\infty}(\overline{\O}^1\times[0,\overline{T}))}.
$$
}

\noindent{\bf \small{Step 3:} Getting contradiction.}
In this step, we have to distinguish two cases:\\

\noindent{\bf \small{Case 1.} ($\overline{t}>0$ and $\overline{s}>0$)}\\
In this case, we have two sub-cases:\\
\noindent{\bf ${\bf(i)}$ If $\overline{x}\in\Omega^1.$} From the definition of $(\overline{x},\overline{t},\overline{s}),$ we see that
$$\overline{u}(\overline{x},t)\leq\varphi^{\overline{u}}(\overline{x},t):=\overline{M}+\overline{v}(\overline{x},\overline{s})+\Psi(\overline{x},t,\overline{s}).$$
Even if $\overline{u}$ has not the required regularity, we see by a simple approximation argument that we can choose
$$\varphi^{\overline{u}}(x,t)=\overline{u}(x,t)\quad \mbox{for}\quad x\neq\overline{x}.$$
Then $\varphi^{\overline{u}}$ is a test function for $\overline{u}$ at $(\overline{x},\overline{t})$
and we deduce the following viscosity inequality
\begin{equation}\label{P_21}
\overline{\b}\frac{\eta}{(\overline{T}-\overline{t})^2}+\overline{\b}\frac{(\overline{t}-\overline{s})}{\delta}\leq\Delta^1[\overline{u}](\overline{x},\overline{t})-\overline{\b}
\lambda\overline{u}(\overline{x},\overline{t}).
\end{equation}
Similarly, we have {\color{black}
$$\overline{v}(\overline{x},s)\geq\varphi^{\overline{v}}(\overline{x},s):=-\overline{M}+\overline{u}(\overline{x},\overline{t})-\Psi(\overline{x},\overline{t},s),$$
}
and we choose
$$\varphi^{\overline{v}}(x,s)=\overline{v}(x,s)\quad \mbox{for}\quad x\neq\overline{x}.$$
We then get the viscosity inequality
\begin{equation}\label{P_22}
\overline{\b}\frac{(\overline{t}-\overline{s})}{\delta}\geq\Delta^1[\overline{v}](\overline{x},\overline{s})-\overline{\b}
\lambda\overline{v}(\overline{x},\overline{s}).
\end{equation}
Subtracting the two viscosity inequalities and setting $w(x)=\overline{u}(x,\overline{t})-\overline{v}(x,\overline{s}),$ we get
\begin{eqnarray*}
0\leq\overline{\b}\frac{\eta}{(\overline{T}-\overline{t})^2}&\leq&\Delta^1[w](\overline{x})-\overline{\beta}
\lambda w(\overline{x})\\
&\leq&\Delta^1[\overline{M}+\Psi(\cdotp,\overline{t},\overline{s})](\overline{x})\\
&\leq&\alpha\sum_{|y-\overline{x}|=1,y\in\mathbb{Z}^n}(|y'|^2-|\overline{x}'|^2)+\gamma\sum_{|y-\overline{x}|=1,y\in\mathbb{Z}^n}(\sqrt{1+y_n}-\sqrt{1+\overline{x}_n})\;{\color{black}=:}\;
A+B,
\end{eqnarray*}
where in the second line {\color{black} we have used the fact that
}$w(y)-\Psi(y,\overline{t},\overline{s})\leq\overline{M}$ for all
$y\in\overline{\O}^1$ with equality for $y=\overline{x}.$ Remark
that, $B=\gamma\zeta(\overline{x}_n)$ with, for $a\geq1$,
$$\zeta(a)=f(a+1)+f(a-1)-2f(a)\qquad\mbox{with}\quad f(a)=\sqrt{1+a},$$
hence we get as in (\ref{eq::comp16}) or (\ref{eq::concave})
{\color{black}
$$\zeta(a)=\sum_{\pm}\int_0^1\,dt\int_0^t\,dsf''(a\pm s)\leq f''(a)$$
} where we have used the concavity of ${\color{black}f''}$. Using
moreover the a priori estimate $(\ref{P_29})$ we have
$\sqrt{1+\overline{x}_n}\leq C/\gamma,$ and we deduce that
{\color{black}
$$B\leq -C'\gamma^4\quad\mbox{with}\quad C'=\frac{1}{4C^3}.$$
}
On the other hand, we have
\begin{equation}\label{estimationonA}
A=\alpha\sum_{|y-\overline{x}|=1,y\in\mathbb{Z}^n}(y'-\overline{x}')(y'+\overline{x}')
\leq2n\alpha(1+2|\overline{x}'|)\leq2n(\alpha+2\sqrt{\alpha}\sqrt{C}),
\end{equation}
where we have used the a priori estimate $(\ref{P_29})$ $(\alpha|\overline{x}'|^2\leq C)$. Finally, if we take $\alpha,\gamma$
small enough so that $2n(\alpha+2\sqrt{\alpha}\sqrt{C})<C'\gamma^4,$ we conclude that $0\leq A+B<0$, and hence a contradiction.\\
\noindent{\bf ${\bf(ii)}$ If $\overline{x}\in\partial\Omega^1.$}
This is a similar case of the latter one where we use the same
arguments. After the same choice of the test function, we arrive at
\begin{equation}\label{P_30}
\frac{\eta}{(\overline{T}-\overline{t})^2}+\frac{\overline{t}-\overline{s}}{\delta}
\leq\overline{F}(\overline{u}(\overline{x},\overline{t}),\overline{t})+D^1[\overline{u}](\overline{x},\overline{t}),
\end{equation}
and
\begin{equation}\label{P_31}
\frac{\overline{t}-\overline{s}}{\delta}\geq\overline{F}(\overline{v}(\overline{x},\overline{s}),\overline{s})+D^1[\overline{v}](\overline{x},\overline{s}).
\end{equation}
Subtracting $(\ref{P_30})$ and
$(\ref{P_31}),$ we infer that

\begin{eqnarray}\label{estimationoneta}
    \frac{\eta}{(\overline{T}-\overline{t})^2}&\leq& {\color{black}\gamma (\sqrt{2}-1)}
+\sum_{|y-\overline{x}|=1,y\in\mathbb{Z}^n {\color{black}, y_n\ge \overline{x}_n}}\alpha(|y'|^2-|\overline{x}'|^2)+[\overline{F}(\overline{u}(\overline{x},\overline{t}),\overline{t})-\overline{F}(\overline{v}(\overline{x},\overline{s}),\overline{s})] \nonumber\\
   &{\color{black}\leq}& {\color{black} \gamma} +\sum_{|y-\overline{x}|=1,y\in\mathbb{Z}^n}\alpha(|y'|^2-|\overline{x}'|^2)+[\overline{F}(\overline{u}(\overline{x},\overline{t}),\overline{t})-\overline{F}(\overline{u}(\overline{x},\overline{t}),\overline{s})] \nonumber\\
   &{}&+\;[\overline{F}(\overline{u}(\overline{x},\overline{t}),\overline{s})-\overline{F}(\overline{v}(\overline{x},\overline{s}),\overline{s})].
\end{eqnarray}
Now, if we take $\lambda\geq
\varepsilon\|F'\|_{L^\infty(\mathbb{R})},$ we get
$\overline{F}'_u\le 0$ and
\begin{equation}\label{estimationonF1}
\overline{F}(\overline{u}(\overline{x},\overline{t}),\overline{s})-\overline{F}(\overline{v}(\overline{x},\overline{s}),\overline{s})\leq0,
\end{equation}
thanks to the fact that
$\overline{u}(\overline{x},\overline{t})-\overline{v}(\overline{x},\overline{s})\geq0.$
Moreover, as $|\overline{F}_t|\leq C_0=
C_0(\eps,\lambda,\|F\|_{W^{1,\infty}(\mathbb{R})}),$ we get
\begin{equation}\label{estimationonF2}
|\overline{F}(\overline{u}(\overline{x},\overline{t}),\overline{t})-\overline{F}(\overline{u}(\overline{x},\overline{t}),\overline{s})|\leq
C_0|\overline{t}-\overline{s}|\leq
C_0\sqrt{2C\delta} ,
\end{equation}
where we have used the a priori estimate $(\ref{P_29})$
$({|\overline{t}-\overline{s}|^2}/{(2\delta)}\leq C).$ Now,
substituting $(\ref{estimationonA}),$ $(\ref{estimationonF1})$ and
$(\ref{estimationonF2})$ into $(\ref{estimationoneta}),$ we infer
that
$$
 0<\frac{\eta}{(T-\overline{t})^2}\leq2n(\alpha+2\sqrt{\alpha}\sqrt{C}) + {\color{black}\gamma}+C_0\sqrt{2C\delta}.
$$
Finally, for $\eta>0$ fixed, we get a contradiction by choosing
$\delta,\alpha$ and $\gamma$ small enough.\\

{\color{black}
\noindent{\bf \small{Case 2.} ($\overline{t}=0$ or $\overline{s}=0$)}\\
For $\eta,\alpha,\gamma$ fixed, we assume that there exists a
sequence $\delta\rightarrow0$ and
$(\overline{x},\overline{t},\overline{s})=(\overline{x}^\delta,\overline{t}^\delta,\overline{s}^\delta)\in\overline{\O}^1\times[0,\overline{T})^2$
such that $\overline{t}^\delta=0$ or $\overline{s}^\delta=0.$ We
deal with the case $\overline{t}^\delta=0$ as the other case
$\overline{s}^\delta=0$ is similar. Then
$$\overline{u}(\overline{x}^\delta,0)-\overline{v}(\overline{x}^\delta,\overline{s}^\delta)=\overline{M}+\Psi(\overline{x}^\delta,0,\overline{s})\geq\frac{M}{4}>0.$$
From $(\ref{P_29}),$ we deduce that, up to an extraction of a
subsequence, we have
$(\overline{x}^\delta,\overline{s}^\delta)\rightarrow(\overline{\overline{x}},0)\in\overline{\O}^1\times[0,\overline{T})$
as $\delta\rightarrow0.$ Therefore for $\beta >0$ we have
$$0<\frac{M}{4}\leq\limsup_{\delta\rightarrow0}\ (\overline{u}(\overline{x}^\delta,0)-\overline{v}(\overline{x}^\delta,\overline{s}^\delta))\leq \overline{u}(\overline{\overline{x}},0)-\overline{v}(\overline{\overline{x}},0)\leq0,$$
where we have used the comparison to the initial condition. This
leads to a contradiction in the case $\beta >0$ or $\beta =0$ and
$\overline{\overline{x}}\in {\color{black}\partial \Omega^1}$. For
the case $\beta=0$ and
$\overline{\overline{x}}\in {\color{black}\Omega^1}$, we get a contradiction exactly as in case 1 i).\\
}

It is worth noticing
that, in the whole proof, we first fix $\eta$, and then we choose
respectively $\gamma$, $\alpha$ and $\delta$ small enough.
$\hfill\square$

\section{Proof of Theorem ${\bf\ref{Comp_pri_con}}$}\label{sec8}

In this section we first present the construction of an auxiliary
function $\xi$ which plays a crucial role in the proof of Theorem
${\ref{Comp_pri_con}}.$ {\color{black}
\begin{lem}${\bf(\mbox{{\bf Auxiliary function}})}$\label{auxiliaryfunction}\\
Let $u\in USC(\overline{\O}\times [0,T)),$ $v\in
LSC(\overline{\O}\times [0,T))$ such that $u\leq C_0,$ $v\geq -C_0,$
$C_0\ge 0,$ and let $M:=\sup_{\overline{\O}\times[0,T)}(u-v)>0.$
Then for all $c>0,$ there exists a
$C^\infty$-function
$\xi_c:\mathbb{R}^n\times[0,T)\rightarrow\mathbb{R}$ with $|\xi_c|\leq 3C_0$ and a positive
constant $a_c>0$ such that for every
$x,y\in\overline{\O}$ and $s,t\in[0,T)$ :
$$v(y,s)\leq\xi_c\left(\frac{x+y}{2},\frac{t+s}{2}\right)\leq
u(x,t),$$
\noindent if $u(x,t)-v(y,s)\geq M-a_c,$ $|x-y|+|t-s|\leq a_c$
and $|x|,|y|\leq c.$
\end{lem}

\noindent {\bf \small{Proof.}}
We essentially revisit the proof of Lemma 5.2 in Barles \cite{Barles}.
Without loss of generality, we can
extend the proof of our result on $[0,T]$ by defining the end point
of $u$ and $v$ as follows:
$$u(x,T)=\limsup_{{\small\begin{aligned}
t<T,\ (y,t)\rightarrow (x,T)
\end{aligned}}}u(y,t),
\quad v(x,T)=\liminf_{{\small\begin{aligned}
t<T,\ (y,t)\rightarrow (x,T)
\end{aligned}}}v(y,t).$$
Let $M:= \sup_{\overline{B}_{2c}\times [0,T]}(u-v)$ with
$\overline{B}_{2c}:=\left\{x\in\overline{\Omega}; |x|\le 2c\right\}$
Let $\mathcal{F}=\{(x,t)\in B_{2c}\times [0,T];\;u(x,t)-v(x,t)=M\}.$
We have $u\le C_0$, $v\ge -C_0$ and $u-v=M$ on $\mathcal{F}$. This
implies that
\begin{equation}\label{eq::borneF}
|u|,|v|\le C_0+M \le 3C_0\quad \mbox{on}\quad \mathcal{F}.
\end{equation}
Moreover it is easy to check that $\mathcal{F}$ is a closed subset
of $\overline{B}_{2c}\times [0,T]$ and that the restriction of $u$
and $v$ to $\mathcal{F}$ are continuous. Therefore, the restriction
to $\mathcal{F}$ of the function $(x,t)\mapsto (u+v)/2$ is also a
continuous function on $\mathcal{F}$ which satisfies
(\ref{eq::borneF}). We may extend this function as a continuous
function in $\mathbb{R}^n\times \R$ (still bounded by $3C_0$) and
then, by standard regularization arguments, there exists a
$C^\infty$-function $\widetilde{\xi}$ such that
\begin{equation}\label{regfunction}
\left|\widetilde{\xi}(x,t)-\frac{u(x,t)+v(x,t)}{2}\right|\leq \frac{M}{8}\quad\mbox{on}\quad\mathcal{F}
\quad \mbox{and}\quad |\tilde{\xi}|\le 3C_0 \quad \mbox{on}\quad \R^n\times\R.
\end{equation}
In order to show the lemma, we now argue by contradiction assuming
that there exist two sequences $(x_{a_c},t_{a_c}),$
$(y_{a_c},s_{a_c})$ such that for $a_c$ small enough,
$u(x_{a_c},t_{a_c})-v(y_{a_c},s_{a_c})\geq M-a_c$ and
$|x_{a_c}-y_{a_c}|+|t_{a_c}-s_{a_c}|\leq a_c$ and such that, say,
$\displaystyle
u(x_{a_c},t_{a_c})-\widetilde{\xi}\left(\frac{x_{a_c}+y_{a_c}}{2},\frac{t_{a_c}+s_{a_c}}{2}\right)<0$.
Extracting, if necessary, subsequences, we may assume without loss
of generality that
$(x_{a_c},t_{a_c}),(y_{a_c},s_{a_c})\rightarrow(\overline{x},\overline{t}).$
Then it is easy to show the convergence of
$u(x_{a_c},t_{a_c})-v(y_{a_c},s_{a_c})$ to
$M=u(\overline{x},\overline{t})-v(\overline{x},\overline{t})$. By
considering the upper semi-continuity of $u$ and the lower
semi-continuity of $v$, this implies, on the one hand,
$u(x_{a_c},t_{a_c})\rightarrow u(\overline{x},\overline{t}),$
$v(y_{a_c},s_{a_c})\rightarrow v(\overline{x},\overline{t})$. On the
other hand, by using the continuity of $\widetilde{\xi}$, we obtain
\begin{equation}\label{contradiction1}
u(\overline{x},\overline{t})-\widetilde{\xi}(\overline{x},\overline{t})=\lim_{a_c\to 0}\left(u(x_{a_c},t_{a_c})-\widetilde{\xi}\left(\frac{x_{a_c}+y_{a_c}}{2},\frac{t_{a_c}+s_{a_c}}{2}\right)\right)\leq0.
\end{equation}
But since $u(\overline{x},\overline{t})-v(\overline{x},\overline{t})=M,$ with $(\overline{x},\overline{t})\in\mathcal{F}$,
we deduce from $(\ref{regfunction})$ that
$$\widetilde{\xi}(\overline{x},\overline{t})
\leq \frac{u(\overline{x},\overline{t})+v(\overline{x},\overline{t})}{2}+\frac{M}{8}
=u(\overline{x},\overline{t}) -\frac{3M}{8}
<u(\overline{x},\overline{t}),$$
which contradicts $(\ref{contradiction1}).$ Finally, we arrive
to the result by taking $\xi_c$ as the restriction of $\widetilde{\xi}$ on $[0,T)$
multiplied by a cut-off function $\psi\in C^\infty(\mathbb{R}^n)$ such that $\psi=1$ on $B_c$ and zero outside $B_{2c}.$ $\hfill\square$\\
}

\noindent{\bf\small{Proof of Theorem $\ref{Comp_pri_con}.$}} The proof is divided into three steps.\\
\noindent{\bf \small{Step 1:} Test function.} In order to replace
the nonlinearity $F$ by a monotone one, we define the new functions
$$\overline{u}:=e^{-\lambda t}u\quad\mbox{and}\quad  \overline{v}:=e^{-\lambda t}v,$$
where $\lambda>0$ is a constant which will be determined later. Obviously, $\overline{u}$ $(\hbox{resp.}\;\overline{v})$ is a sub-solution (resp. super-solution) of the following problem \begin{equation}\label{newPN}\left\{
  \begin{array}{ll}
   \,\,\displaystyle{\beta\overline{u}_t=\Delta[\overline{u}]-\beta\lambda\overline{u}}&\displaystyle{\quad\hbox{in}\quad\Omega\times(0,T)},\\
   {}\\
   \displaystyle{\overline{u}_t=\frac{\partial\overline{u}}{\partial x_n}+\overline{F}(\overline{u},t)}&\displaystyle{\quad\hbox{on}\quad\partial\Omega\times(0,T)},
  \end{array}
\right.
\end{equation}
and
\begin{equation}\label{newInit_cond_PN}
\overline{u}(x,0)=u_{0}(x)\quad \mbox{for} \quad
\left\{
\begin{aligned}
&x\in \overline{\O}\quad
&\mbox{if}& \quad \b> 0\\
& x\in \partial \O\quad &\mbox{if}& \quad \b=0,
\end{aligned}
\right.
\end{equation}
with
$$\overline{F}(\overline{u},t)= e^{-\lambda t}F(e^{\lambda t}\overline{u})-\lambda\overline{u}.$$
Let us assume that
$M:=\sup_{\overline{\Omega}\times [0,T)}(\overline{u}-\overline{v})>0$
and let us exhibit a contradiction. By the definition of the
supremum, there exists
$P^{k_0}=(x^{k_0},t^{k_0})\in\overline{\O}\times[0,T),$ for some
index $k_0,$ such that
\begin{equation}\label{levelk_0}
\overline{u}(P^{k_0})-\overline{v}(P^{k_0})\geq \frac{M}{2}.
\end{equation}

{\color{black} Let us introduce the following constant:
$$C_*=\|\overline{u}\|_{L^{\infty}(\overline{\O}\times[0,T))}+\|\overline{v}\|_{L^{\infty}(\overline{\O}\times[0,T))}.$$
We now take the following notation: $x=(x',x_n)$ with
$x'=(x_1,\dots,x_{n-1})$, and for $\alpha,\gamma,\eta >0$ (to be
fixed later), we can approximate the functions $\overline{u}$,
$\overline{v}$ by the functions:
\begin{equation}\label{P_180}
\tilde{u}(x,t):=\overline{u}(x,t)-\alpha|x'|^2 -\gamma\sqrt{1+x_n}-\frac{\eta}{T-t} \le C_*,
\end{equation}
\begin{equation}\label{P_1800}
\tilde{v}(y,s):=\overline{v}(y,s)+\alpha|y'|^2 +\gamma\sqrt{1+y_n}+\frac{\eta}{T-s} \ge -C_*
\end{equation}
and define
$$\widetilde{M}:=\sup_{x\in\overline{\O}, t\in[0,T)}(\tilde{u}(x,t)-\tilde{v}(x,t)).$$
In order to dedouble the variables in space and time, following the proof of Theorem 3.2 in \cite{Barles},
we define the function $\tilde{\Phi}:(\overline{\O}\times[0,T))^2\rightarrow\mathbb{R}$ by:
\begin{eqnarray*}
  \tilde{\Phi}(x,t,y,s)&:=&\frac{(x-y)^2}{\varepsilon^2}+\frac{B(x_n-y_n)^2}{\varepsilon^2}+\frac{(t-s)^2}{\delta^2}+\frac{2(t-s)(x_n-y_n)}{\delta^2}\\
 &{}&-(x_n-y_n)\overline{F}\left(\xi\left(\frac{x+y}{2},\frac{t+s}{2}\right),\frac{t+s}{2}\right),
\end{eqnarray*}
where the parameters $B,\eps,\delta>0$ will be chosen later.
Moreover $c>0$ is a constant that will be defined later
(only depending on $\alpha,\gamma,\lambda,C_*,\|F\|_{L^\infty(\mathbb{R})}$) and the smooth function
\begin{equation}\label{eq::xic}
\xi=\xi_c
\end{equation}
is  the auxiliary function associated to $\tilde{u}$ and $\tilde{v}$, and
given by Lemma $\ref{auxiliaryfunction}$, which shows in particular the following estimate
\begin{equation}\label{eq::bornexi}
|\xi|\le 3 C_*
\end{equation}

We note that $ \tilde{\Phi}(x,t,x,t)=0.$
Then we set
\begin{equation}\label{P_1801}
\overline{M}=\overline{M}_{B,\eps,\delta}:=\sup_{x,y\in\overline{\O},\  t,s\in[0,T)}(\tilde{u}(x,t)-\tilde{v}(y,s)-\tilde{\Phi}(x,t,y,s))\ge \tilde{M}.
\end{equation}
}

\noindent{\bf \small{Step 2:} A priori estimates.} By choosing
$\eta,\alpha,\gamma$ small enough such that
$$\frac{\eta}{T-t^{k_0}}+\alpha|(x^{k_0})'|^2+\gamma\sqrt{1+x^{k_0}_n}\leq\frac{M}{8},$$
it follows from (\ref{levelk_0}) that
\begin{equation}\label{estimationonM}
\overline{M}\geq\widetilde{M}\geq \tilde{u}(P^{k_0})-\tilde{v}(P^{k_0})\geq\frac{M}{4}>0.
\end{equation}
Hence, from the definition of $\overline{M},$ there exist sequences
$\overline{x}^k,\overline{y}^k\in\overline{\O},$
$\overline{t}^k,\overline{s}^k\in[0,T)$ such that
$$\tilde{u}(\overline{x}^k,\overline{t}^k)-\tilde{v}(\overline{y}^k,\overline{s}^k)-\tilde{\Phi}(\overline{x}^k,\overline{t}^k,\overline{y}^k,\overline{s}^k)\longrightarrow\overline{M}>0\quad\mbox{as}\quad k\rightarrow\infty,$$
and, by taking $k$ large enough, we deduce that
\begin{equation}\label{P_190}
\tilde{u}(\overline{x}^k,\overline{t}^k)-\tilde{v}(\overline{y}^k,\overline{s}^k)-\tilde{\Phi}(\overline{x}^k,\overline{t}^k,\overline{y}^k,\overline{s}^k)\geq\frac{\overline{M}}{2}>0.
\end{equation}
{\color{black} From (\ref{eq::bornexi}), it follows that
$$|\overline{F}(\xi(\frac{\overline{x}^k+\overline{y}^k}{2},\frac{\overline{t}^k+\overline{s}^k}{2}),\frac{\overline{t}^k+\overline{s}^k}{2})|\leq C_1=C_1(\lambda,C_*,\|F\|_{L^\infty(\mathbb{R})}).$$
Then if we take $\eps\le \delta \le 1$, Young's inequality yields
\begin{equation}\label{Young1}
|(\overline{x}^k_n-\overline{y}^k_n)\overline{F}(\xi(\frac{\overline{x}^k+\overline{y}^k}{2},\frac{\overline{t}^k+\overline{s}^k}{2}),\frac{\overline{t}^k+\overline{s}^k}{2})|\leq \frac{(\overline{x}^k_n-\overline{y}^k_n)^2}{\varepsilon^2}+\frac{C_1^2\varepsilon^2}{4}\leq \frac{(\overline{x}^k_n-\overline{y}^k_n)^2}{\varepsilon^2}+\frac{C_1^2\delta^2}{4},
\end{equation}
and
\begin{equation}\label{Young2}
\frac{|2(\overline{t}^k-\overline{s}^k)(\overline{x}_n^k-\overline{y}_n^k)|}{\delta^2}\leq \frac{1}{2}\frac{(\overline{t}^k-\overline{s}^k)^2}{\delta^2}+2\frac{(\overline{x}^k_n-\overline{y}^k_n)^2}{\varepsilon^2}.
\end{equation}
Using  (\ref{Young1}), (\ref{Young2}), we deduce that for $B\ge 3$
$$\tilde{\Phi}(\overline{x}^k,\overline{t}^k,\overline{y}^k,\overline{s}^k) \ge \frac{(x-y)^2}{\varepsilon^2}+\frac{(B-3)(x_n-y_n)^2}{\varepsilon^2}+\frac{(t-s)^2}{2\delta^2}-\frac{C_1^2\delta^2}{4} \ge -\frac{C_1^2\delta^2}{4}.$$
Using (\ref{P_190}), we conclude that
$$\frac{(\overline{x}^k-\overline{y}^k)^2}{\varepsilon^2}+\frac{(\overline{t}^k-\overline{s}^k)^2}{2\delta^2}+\alpha(|(\overline{x}^k)'|^2+|(\overline{y}^k)'|^2)+\gamma\sqrt{1+\overline{x}^k_n}+\gamma\sqrt{1+\overline{y}^k_n}+\frac{\eta}{T-\overline{t}^k}+\frac{\eta}{T-\overline{s}^k}\leq C_{**},$$
with (for $\delta\le 1$)
$$C_{**}=C_* + \frac{C_1^2}{4} = C_{**}(\lambda, C_*,\|F\|_{L^\infty(\mathbb{R})}).$$
Then, up to an extraction of
a subsequence, we have
$$(\overline{x}^k,\overline{t}^k,\overline{y}^k,\overline{s}^k)\rightarrow(\overline{x},\overline{t},\overline{y},\overline{s})\in(\overline{\O})^2\times[0,T)^2\quad\mbox{as}\;\;k\rightarrow\infty,$$
with
\begin{equation}\label{P_280}
\tilde{u}(\overline{x},\overline{t})-\tilde{v}(\overline{y},\overline{s})-\tilde{\Phi}(\overline{x},\overline{t},\overline{y},\overline{s})=\overline{M}>0 \quad \mbox{and}\quad \tilde{\Phi}(\overline{x},\overline{t},\overline{y},\overline{s})\ge
-\frac{C_1^2\delta^2}{4}
\end{equation}
and the same a priori estimate
\begin{equation}\label{P_290}
\frac{(\overline{x}-\overline{y})^2}{\varepsilon^2}+\frac{(\overline{t}-\overline{s})^2}{2\delta^2}+\alpha(|\overline{x}'|^2+|\overline{y}'|^2)+\gamma\sqrt{1+\overline{x}_n}+\gamma\sqrt{1+\overline{y}_n}+\frac{\eta}{T-\overline{t}}+\frac{\eta}{T-\overline{s}}\leq {C}_{**}.
\end{equation}
}

\noindent{\bf \small{Step 3:} Getting contradiction.} In order to get a contradiction, we have to distinguish two cases:\\

\noindent{\bf \small{Case 1.} ($\overline{t}>0$ and $\overline{s}>0$)}\\
\noindent{\bf ${\bf(i)}$ If $\overline{x}\in\partial\Omega$ or
$\overline{y}\in\partial\Omega$.} We only study the case
$\overline{x}\in\partial\Omega$ as it is similar for
$\overline{y}\in\partial\Omega$. Let
$\varphi^{\overline{u}}:\overline{\O}\times[0,T)\rightarrow\mathbb{R}$
be the function defined by {\color{black}
$$\varphi^{\overline{u}}(x,t):=
\displaystyle{\overline{M}+\tilde{v}(\overline{y},\overline{s})
+\tilde{\Phi}({x},t,\overline{y},\overline{s})+\alpha|{x}'|^2
+\gamma\sqrt{1+{x}_n}+\frac{\eta}{T-t}}.$$ } Using $(\ref{P_280}),$
we see that $\varphi^{\overline{u}}$ is a test function for
$\overline{u}$ and hence we obtain
$$ \min \left\{\beta \varphi^{\overline{u}}_{t}(P_{0})-\Delta \varphi^{\overline{u}}(P_{0})-\beta\lambda\overline{u} \;,\;
\varphi^{\overline{u}}_{t}(P_0) - \overline{F}(\varphi^{\overline{u}}(P_0),\overline{t}) - \frac{\partial
  \varphi^{\overline{u}}}{\partial x_{n}}(P_0)\right\}\leq 0,$$
 where $P_0=(\overline{x},\overline{t}).$ We note that the case $\beta \varphi^{\overline{u}}_{t}(P_{0})-\Delta \varphi^{\overline{u}}(P_{0})-\beta\lambda\overline{u} \leq 0$ can be treated as in the sub-case ${\bf(ii)}$ below, so it is sufficient to study the case where $\varphi^{\overline{u}}_{t}(P_0) - \overline{F}(\varphi^{\overline{u}}(P_0),\overline{t}) - \frac{\partial
  \varphi^{\overline{u}}}{\partial x_{n}}(P_0)\leq0,$ which implies, with $\overline{x}_n=0$:
{\color{black}
\begin{eqnarray}\label{inequalityA}
&{}&\frac{\eta}{(T-\overline{t})^2}+\frac{2(\overline{x}_n-\overline{y}_n)}{\delta^2}
-\frac{\gamma}{2\sqrt{1+\overline{x}_n}}
+\;\overline{F}\left(\xi\left(\frac{\overline{x}+\overline{y}}{2},\frac{\overline{t}+\overline{s}}{2}\right),\frac{\overline{t}+\overline{s}}{2}\right)-\overline{F}(\overline{u}(\overline{x},\overline{t}),\overline{t})\nonumber\\
&{}&
-\frac{2(1+B)(\overline{x}_n-\overline{y}_n)}{\varepsilon^2}
-\frac{(\overline{x}_n-\overline{y}_n)}{2}(\overline{F})^{'}_u.\xi^{'}_t
+\frac{(\overline{x}_n-\overline{y}_n)}{2}(\overline{F})^{'}_u.\xi^{'}_x
-\frac{(\overline{x}_n-\overline{y}_n)}{2}(\overline{F})^{'}_t\leq0,
\end{eqnarray}
where
\begin{equation}\label{bounded1}
\xi^{'}_t=\xi^{'}_t\left(\frac{\overline{x}+\overline{y}}{2},\frac{\overline{t}+\overline{s}}{2}\right),\qquad(\overline{F})^{'}_t=(\overline{F})^{'}_t\left(\xi\left(\frac{\overline{x}+\overline{y}}{2},\frac{\overline{t}+\overline{s}}{2}\right),\frac{\overline{t}+\overline{s}}{2}\right)
\end{equation}
and
\begin{equation}\label{bounded2}
\xi^{'}_x=\xi^{'}_x\left(\frac{\overline{x}+\overline{y}}{2},\frac{\overline{t}+\overline{s}}{2}\right),\qquad (\overline{F})^{'}_u=(\overline{F})^{'}_u\left(\xi\left(\frac{\overline{x}+\overline{y}}{2},\frac{\overline{t}+\overline{s}}{2}\right),\frac{\overline{t}+\overline{s}}{2}\right)
\end{equation}
are the first partial derivatives of $\xi$ and $\overline{F}.$ On
the one hand, from {\color{black}(\ref{P_280}) and (\ref{P_1801})},
we deduce that
\begin{equation}\label{eq::u-v}
\tilde{u}(\overline{x},\overline{t})-\tilde{v}(\overline{y},\overline{s}) \ge \tilde{M}-\frac{C_1^2\delta^2}{4}.
\end{equation}
From (\ref{P_290}) we deduce that
\begin{equation}\label{eq::xybound}
|\overline{x}|, |\overline{y}| \le C_2=C_2(\alpha,\gamma,C_{**})
\quad \mbox{and}\quad |\overline{x}-\overline{y}| \le \varepsilon \sqrt{C_{**}},
\quad  |\overline{t}-\overline{s}| \le \delta \sqrt{2C_{**}}.
\end{equation}
Then in (\ref{eq::xic}), we can choose
$$c=c(\alpha,\gamma, \lambda, C_*, \|F\|_{L^\infty(\mathbb{R})}):=C_2$$
and Lemma \ref{auxiliaryfunction} gives the existence of a number $a_c >0$.
Then choosing $0<\varepsilon \le \delta$ small enough, we deduce from (\ref{eq::u-v}) and (\ref{eq::xybound}) that
$$\tilde{u}(\overline{x},\overline{t})-\tilde{v}(\overline{y},\overline{s}) \ge \tilde{M}-a_c,\quad
|\overline{x}-\overline{y}| +|\overline{t}-\overline{s}| \le a_c
\quad \mbox{and}\quad |\overline{x}|,|\overline{y}| \le c.$$
Therefore, we deduce from Lemma \ref{auxiliaryfunction} that
$$\xi\left(\frac{\overline{x}+\overline{y}}{2},\frac{\overline{t}+\overline{s}}{2}\right)
\le \tilde{u}(\overline{x},\overline{t}) \le \overline{u}(\overline{x},\overline{t}).$$
On the other hand, for the choice $\lambda\ge \|F'\|_{L^\infty(\mathbb{R})})$, we have $\overline{F}'_{u}\le 0$
and then we can estimate
\begin{equation}\label{inequalityB}
\overline{F}\left(\xi\left(\frac{\overline{x}+\overline{y}}{2},\frac{\overline{t}+\overline{s}}{2}\right),\frac{\overline{t}+\overline{s}}{2}\right)-\overline{F}(\overline{u}(\overline{x},\overline{t}),\overline{t})
\geq -C|\overline{t}-\overline{s}|.
\end{equation}
Now, as the terms in $(\ref{bounded1})$ and $(\ref{bounded2})$ are
bounded, we conclude from {\color{black} $(\ref{inequalityA}),$
$(\ref{inequalityB})$ and $(\ref{P_290})$} that
$$\frac{\eta}{T^2}-\frac{\gamma}{2}\leq O(\frac{\varepsilon}{\delta^2})+ O(\varepsilon)+O(\delta).$$
Moreover, for
$\eta>0$ fixed, we get the contradiction if we take
$\eps,\delta,\gamma$ small enough so that
$\eps=\eps(\delta)\le\delta^3\le 1$ and ${\gamma}<{\eta}/{T^2}$. We
note that in the case of
$(\overline{x},\overline{y})\in\O\times\partial\O,$ the test
function for $\overline{v}$ is given by
$$\varphi^{\overline{v}}(y,s):=
\displaystyle{-\overline{M}+\tilde{u}(\overline{x},\overline{t})-\tilde{\Phi}(\overline{x},\overline{t},{y},s)-\alpha|{y}'|^2
-\gamma\sqrt{1+{y}_n}-\frac{\eta}{T-s}}.$$ }

{\color{black}
\noindent{\bf ${\bf(ii)}$ If $\overline{x}\in\Omega$ and $\overline{y}\in\Omega$.}
We know from $(\ref{P_280})$ that
$\overline{u}(\overline{x},\overline{t})-\overline{v}(\overline{y},\overline{s})
-\Phi(\overline{x},\overline{t},\overline{y},\overline{s})$
has a local maximum at $(\overline{x},\overline{y},\overline{t},\overline{s}),$ where
$$\Phi(\overline{x},\overline{t},\overline{y},\overline{s})
:=\tilde{\Phi}(\overline{x},\overline{t},\overline{y},\overline{s})
+\alpha(|\overline{x}'|^2+|\overline{y}'|^2)
+\gamma\left(\sqrt{1+\overline{x}_n}+\sqrt{1+\overline{y}_n}\right)
+\frac{\eta}{T-\overline{t}}+\frac{\eta}{T-\overline{s}}.$$

 Then it is natural to apply the classical Ishii's Lemma in the elliptic case
with the new coordinates $(\tilde{x}=(x,t),\tilde{y}=(y,s))$ and this is what we do.
Indeed, we only use a corollary of Ishii's Lemma,
namely Corollary \ref{PVI} which is given in the Appendix.
Applying Corollary \ref{PVI}, we get for every $\mu>0$ satisfying $\mu
\overline{{A}}<\overline{I},$ (with $\overline{{A}}$ defined in (\ref{eq::Abar}) and $\overline{I}$ is the identity matrix of $\R^{2(n+1)}$),
the existence of symmetric $n\times n$ matrices $X,Y$ such that
\begin{equation}\label{estimation8}
(\tau_1,p_1,
\left(\begin{array}{cc}
X & * \\
* & * \\
\end{array}\right)
)\in\overline{\mathcal{D}}^+\overline{u}(\overline{x},\overline{t}),
\;\;(\tau_2,p_2,
\left(\begin{array}{cc}
Y & * \\
* & * \\
\end{array}\right)
)\in\overline{\mathcal{D}}^-\overline{v}(\overline{y},\overline{s}),
\end{equation}
and
\begin{equation}\label{estimation3}
 -(\frac{1}{\mu}+\|\overline{{A}}\|)\hat{I}\leq
\left(\begin{array}{cc}
X & 0 \\
0 & -Y \\
\end{array}\right)
\leq {A}+2\mu({A}^2+A_1\cdot A_1^T),
\end{equation}
where $\hat{I}$ is the identity matrix of $\R^{2n}$ and
where all the quantities are computed at the point $(\overline{x},\overline{t},\overline{y},\overline{s})$:
$\tau_1:=\partial_t\Phi,$
$p_1:=D_x\Phi,$ $\tau_2:=-\partial_s\Phi,$
$p_2:=-D_y\Phi,$
and
$${A}=
\left(\begin{array}{cc}
D^2_{xx}\Phi & D^2_{xy}\Phi \\
D^2_{yx}\Phi & D^2_{yy}\Phi \\
\end{array}\right)
\quad \mbox{and}\quad
A_1=
\left(\begin{array}{cc}
D^2_{xt}\Phi & D^2_{xs}\Phi \\
D^2_{yt}\Phi & D^2_{ys}\Phi \\
\end{array}\right).$$
A simple computation shows that
\begin{eqnarray*}
 {\color{black} A} &=& \frac{2}{\varepsilon^2}\left(
                                            \begin{array}{c|c}
                                            I & -I \\
                                            \hline
                                              -I & I \\
                                            \end{array}
                                          \right)+\frac{2B}{\varepsilon^2}\left(
                                            \begin{array}{c|c}
                                            I_n
                                            &  -I_n \\
                                            \hline
                                               -I_n &  I_n\\
                                            \end{array}
                                            \right)+ 2\alpha\left(
                                            \begin{array}{c|c}
                                            I'
                                            &  0 \\
                                            \hline
                                              0 &  I'\\
                                            \end{array}
                                            \right)-\frac{\gamma}{4}\left(
                                            \begin{array}{c|c}
                                            \frac{1}{(1+\overline{x}_n)^{3/2}}I_n
                                            &  0 \\
                                            \hline
                                              0 &  \frac{1}{(1+\overline{y}_n)^{3/2}}I_n\\
                                            \end{array}
                                            \right)\\
   &-& {\tiny{ \left(
                                                      \begin{array}{c|c}
                                                        \left(
                                                          \begin{array}{ccc|c}
                                                            0 &  &     \\
                                                             &  \cdotp&  &P'  \\
                                                             &  & 0 &  \\
                                                             \hline
                                                             &  P' &  & 2P_n \\
                                                          \end{array}
                                                        \right)
                                                         & \left(
                                                          \begin{array}{ccc|c}
                                                            0 &  &  \\
                                                             & \cdotp &  &-P'  \\
                                                           &  & 0 &  \\
                                                             \hline
                                                             &   P' &  & 0 \\
                                                          \end{array}
                                                        \right) \\
                                                         \hline
                                                        \left(
                                                          \begin{array}{ccc|c}
                                                            0 &   &  \\
                                                             & \cdotp&  &P'  \\
                                                              &  & 0 &  \\
                                                             \hline
                                                             &  -P' &  & 0 \\
                                                          \end{array}
                                                        \right) & \left(
                                                          \begin{array}{ccc|c}
                                                            0 &   &  \\
                                                             & \cdotp &   &-P' \\
                                                            &  & 0 &  \\
                                                             \hline
                                                             &  -P' &  & -2P_n \\
                                                          \end{array}
                                                        \right) \\
                                                      \end{array}
                                                    \right)}}
                                              -(\overline{x}_n-\overline{y}_n)\left(
                                       \begin{array}{c|c}
                                         D^2\tilde{F} & D^2\tilde{F} \\
                                         \hline
                                         D^2\tilde{F}& D^2\tilde{F}\\
                                       \end{array}
                                     \right).
\end{eqnarray*}
Here $I$ is the identity matrix of $\R^n$ and for all $i,j \in
\left\{1,...,n\right\}$, $(I_n)_{i,j}=1$ if $i=j=n$ and $0$
otherwise, and $(I')_{i,j}:=1$ if $i=j\leq n-1$ and $0$ otherwise.
Moreover
$$\tilde{F}(x,y):=\overline{F}(\xi(\frac{\overline{t}+\overline{s}}{2},\frac{{x}+{y}}{2}),
\frac{\overline{t}+\overline{s}}{2}),\quad
P':=D_{x'}\tilde{F}=D_{y'}\tilde{F},\quad
P_n:=D_{x_n}\tilde{F}=D_{y_n}\tilde{F},$$
for all $x=(x',x_n)\in\mathbb{R}^{n-1}\times\mathbb{R}$ and $y=(y',y_n)\in\mathbb{R}^{n-1}\times\mathbb{R}.$
Then we can write the viscosity inequalities for the limit sub/superdifferentials
$\overline{\mathcal{D}}^-\overline{v}(\overline{y},\overline{s})$
and $\overline{\mathcal{D}}^+\overline{u}(\overline{x},\overline{t})$. This gives
$$\beta\tau_1+\beta\lambda \overline{u}(\overline{x},\overline{t})-\mbox{tr}(X)\leq0$$
and
$$\beta\tau_2+\beta\lambda \overline{v}(\overline{y},\overline{s})-\mbox{tr}(Y)\geq0,$$
where $\mbox{tr}$ is the trace of the appropriate matrix.
Substracting these viscosity inequalities and using
$(\ref{P_280}),$ we get
$$\beta(\partial_t\Phi+\partial_s\Phi)+\beta\lambda(\Phi+\overline{M})\leq
\mbox{tr}(X-Y),$$
i.e.
$$\beta\left(\frac{\eta}{(T-\overline{t})^2}+\frac{\eta}{(T-\overline{s})^2}-(\overline{x}_n-\overline{y}_n)((\overline{F})'_t+(\overline{F})'_u\xi'_t)\right)+\beta\lambda(\Phi+\overline{M})\leq
\mbox{tr}(X-Y).$$
Moreover, using the
definition of $\Phi,$ we conclude that
\begin{eqnarray}\label{inequality1}
&{}&\beta\left(\frac{\eta}{(T-\overline{t})^2}+\frac{\eta}{(T-\overline{s})^2}
-(\overline{x}_n-\overline{y}_n)((\overline{F})'_t+(\overline{F})'_u\xi'_t)\right)
+\beta\lambda(\overline{M}-(\overline{x}_n-\overline{y}_n)\overline{F})\nonumber\\
&{}&+\;\beta\lambda\left(\frac{B(\overline{x}_n-\overline{y}_n)^2}{\varepsilon^2}+\frac{(\overline{t}-\overline{s})^2}{\delta^2}+\frac{2(\overline{t}-\overline{s})(\overline{x}_n-\overline{y}_n)}{\delta^2}\right)\nonumber\\
&{}&\leq \mbox{tr}(X-Y).
\end{eqnarray}
Using the fact that $B\ge 1$ and $\overline{M}>0$,  we deduce that
$$2\beta\frac{\eta}{T^2}-\beta(\overline{x}_n-\overline{y}_n)\left((\overline{F})'_t+(\overline{F})'_u\xi'_t+\lambda
\overline{F}\right)\leq \mbox{tr}(X-Y).$$
From (\ref{eq::xybound}), we deduce that
\begin{equation}\label{eq::tr0}
2\beta\frac{\eta}{T^2}\leq
\mbox{tr}(X-Y)+O(\varepsilon).
\end{equation}
On the other hand, taking the matrix inequality (\ref{estimation3}) on the vectors
$\left(\begin{array}{c}
\xi \\
\xi
\end{array}\right)^T$ and $\left(\begin{array}{c}
\xi \\
\xi
\end{array}\right)$, we get
$$\xi ^T (X-Y)\xi \le K_1 + K_2$$
with
$$K_1= \left(\begin{array}{c}
\xi \\
\xi
\end{array}\right)^T A \left(\begin{array}{c}
\xi \\
\xi
\end{array}\right) = 4\alpha |\xi'|^2 -\frac{\gamma}{4}\left(
\frac{1}{(1+\overline{x}_n)^{\frac32}}+\frac{1}{(1+\overline{y}_n)^{\frac32}}\right)\xi_n^2
- 4(\overline{x}_n-\overline{y}_n)\xi^T D^2\tilde{F}\xi$$
and
$$K_2= 2\mu \left(\begin{array}{c}
\xi \\
\xi
\end{array}\right)^T (A^2 + A_1\cdot A_1^T) \left(\begin{array}{c}
\xi \\
\xi
\end{array}\right)\le 4\mu {\color{black}\|}A^2 + A_1\cdot A_1^T{\color{black}\|}\  |\xi|^2.$$
Using again the a priori estimate $(\ref{P_290}),$ we have
$$-\frac{\gamma}{4}\left\{\frac{1}{(1+\overline{x}_n)^{3/2}}+\frac{1}{(1+\overline{y}_n)^{3/2}}\right\}\leq
-C'\gamma^4,$$
and we conclude that
\begin{equation}\label{inequality5}
\mbox{tr}(X-Y)\leq
4\alpha(n-1)-C'\gamma^4-4(\overline{x}_n-\overline{y}_n)\mbox{tr}(D^2\tilde{F})+4\mu
n {\color{black}\|}A^2 + A_1\cdot A_1^T {\color{black}\|}.
\end{equation}
From (\ref{eq::tr0}) and the bounds on $F$ in $W^{2,\infty}(\R)$, we deduce that
$$2\beta\frac{\eta}{T^2}\le 4\alpha(n-1)-C'\gamma^4 +O(\varepsilon) +4\mu n {\color{black}\|}A^2 + A_1\cdot A_1^T {\color{black}\|}.$$
Because all the terms are independent on $\mu$, we can first take the limit $\mu\to 0$.
Recall that $\beta\ge 0$, and then we can not use $\eta$ to get a contradiction as usual.
We then have to get a contradiction only using the following inequality:
$$0\le  4\alpha(n-1)-C'\gamma^4 +O(\varepsilon).$$
To this end, for given $\gamma>0$, we get a contradiction choosing
$\alpha, \eps$ small enough such that $\eps\ll \gamma^4$ and
$8\alpha(n-1)<C^{'}\gamma^4.$\\
}

{\color{black}
\noindent{\bf \small{Case 2.} ($\overline{t}=0$ or $\overline{s}=0$)}\\
By fixing the parameters $\gamma$, $\eta$ and $\alpha$, we assume
that there exists a sequence $\delta\rightarrow0$ and
$(\overline{x},\overline{t},\overline{y},\overline{s})
=(\overline{x}^\delta,\overline{t}^\delta,\overline{y}^\delta,\overline{s}^\delta)\in
(\overline{\Omega}\times [0,T))^2$ such that $\overline{s}^\delta=0$
or $\overline{t}^\delta=0.$ It is sufficient to study the case
$\overline{t}^\delta=0$ as the other case can be deduced similarly.
Thus
$$\tilde{u}(\overline{x}^\delta,0)-\tilde{v}(\overline{y}^\delta,\overline{s}^\delta)=\overline{M}+\tilde{\Phi}(\overline{x}^\delta,0,\overline{y}^\delta,\overline{s}^\delta),$$
and we deduce from $(\ref{P_290})$ (with $\varepsilon\le \delta^3\le
1$) that, up to extraction of a subsequence, we have
$(\overline{x}^\delta,0,\overline{y}^\delta,\overline{s}^\delta)
\rightarrow(\overline{\overline{x}},0,\overline{\overline{x}},0)\in
(\overline{\Omega}\times [0,T))^2$ as $\delta\rightarrow0$. Hence in
the case $\beta>0,$ we conclude that (using (\ref{P_280}))
$$0<\overline{M}
\leq\limsup_{\delta\rightarrow0}\ (\overline{M}+\tilde{\Phi}(\overline{x}^\delta,0,\overline{y}^\delta,\overline{s}^\delta))
=\limsup_{\delta\rightarrow0}\ (\tilde{u}(\overline{x}^\delta,0)-\tilde{v}(\overline{y}^\delta,\overline{s}^\delta))
\leq\overline{u}(\overline{\overline{x}},0)-\overline{v}(\overline{\overline{x}},0)\leq0,$$
where we have used the comparison to the initial condition.
This gives a contradiction in the case $\beta>0$ and in the case $\beta=0$
if $\overline{\overline{x}}\in \partial \Omega$.
In the case $\beta=0$ and $\overline{\overline{x}}\in \Omega$, we get a contradiction exactly as in the case 1 i).\\
}

In the whole proof, we first fix
$\eta$ and $B$. After that, we fix respectively  $\gamma$, $\alpha$,
$\delta$, $\eps$ (and $\mu$). $\hfill{\Box}$

\section{Appendix}

This appendix is dedicated to the proof of a
corollary of Ishii's lemma (Corollary \ref{PVI}), which is used in the proof of Theorem \ref{Comp_pri_con}.
First, we recall the elliptic
sub and superdifferentials of semi-continuous functions and the classical Ishii's Lemma.
In all what follows, we denote by $\S^n$ the set of symmetric
$n\times n$ matrices.

\begin{definition}${\bf(\mbox{{\bf Elliptic sub and superdifferential of order two}})}$\label{subsuperelliptic}\\
Let $U$ be a locally compact subset of $\mathbb{R}^n$ and $u\in
USC(U).$ Then the superdifferential $\mathcal{D}^{+}u$ of order
two of the function $u$ is defined by: $(p,X)\in \mathbb{R}^n\times
\S^n$ belongs to $\mathcal{D}^{+}u(x)$ if $x\in U$ and
$$u(y)\leq u(x)+\langle p,y-x\rangle+\frac{1}{2}\langle X(y-x),y-x\rangle+o(|y-x|^2)$$
as $U\ni y\rightarrow x.$ In a similar way, we define the
subdifferential of order two by
$\mathcal{D}^{-}u=-\mathcal{D}^{+}(-u).$ We also define:
$$\overline{\mathcal{D}}^{+}u(x)=\left\{
\begin{array}{ll}
\,\,\displaystyle{(p,X)\in\mathbb{R}^n\times
\S^n,\;\exists(x_n,p_n,X_n)\in U\times\mathbb{R}^n\times
\S^n}\\
\displaystyle{\hbox{such
that}\;(p_n,X_n)\in\mathcal{D}^{+}u(x_n)}\\
\displaystyle{\hbox{and}\;(x_n,u(x_n),p_n,X_n)\rightarrow(x,u(x),p,X)}
\end{array}
\right\}.$$ The set $\overline{\mathcal{D}}^{-}u(x)$ is defined in
a similar way.
\end{definition}

Now recall the classical elliptic version of Ishii's Lemma.

\begin{lemma}\label{elliptic}${\bf(\mbox{{\bf Elliptic version of Ishii's lemma}})}$\\
Let $U$ and $V$ be a locally compact subsets of $\mathbb{R}^n,$
$u\in USC(U)$ and $v\in LSC(V).$ Let $\varphi:U\times
V\rightarrow\mathbb{R}$ be of class $C^2.$ Assume that $(x,y)\mapsto
u(x)-v(y)-\varphi(x,y)$ reaches a local maximum at
$(\overline{x},\overline{y})\in U\times V.$ We note
$p_1=D_x\varphi(\overline{x},\overline{y}),$
$p_2=-D_y\varphi(\overline{x},\overline{y})$ and
$A=D^2\varphi(\overline{x},\overline{y}).$ Then, for every $\mu>0$
such that $\mu A<\hat{I},$ there exists $X,Y\in \S^n$ such that:
\begin{equation}\label{ishii1}
\begin{array}{l}
 (p_1,X)\in\overline{\mathcal{D}}^{+}u(\overline{x}),\;\;(p_2,Y)\in\overline{\mathcal{D}}^{-}v(\overline{y}),\\
   \displaystyle -(\frac{1}{\mu}+\|A\|)\hat{I}\leq \left(
                                                      \begin{array}{cc}
                                                        X & 0 \\
                                                        0 & -Y \\
                                                      \end{array}
                                                    \right)\leq
                                                    A+\mu A^2,
\end{array}\end{equation}
where $\hat{I}$ is the identity matrix of $\R^{2n}$.
The norm of the symmetric matrix $A$ used in $(\ref{ishii1})$ is
$$\|A\|=\sup\{|\lambda|:\;\lambda\;\hbox{is an eigenvalue
of}\;A\}=\sup\{|<A\xi,\xi>|:|\xi|\leq1\}.$$
\end{lemma}
\noindent For the proof, we refer the reader to Theorem $3.2$ in the
User's Guide \cite{CIL}.

Then we have:

{\color{black}
\begin{corollary}\label{PVI}${\bf(\mbox{{\bf Consequence of Ishii's lemma}})}$\\
Given $T>0.$ Let $U$ and $V$ be a locally compact subsets of
$\mathbb{R}^n,$ $u\in USC(U\times[0,T))$ and $v\in
LSC(V\times[0,T)).$ Let $\Phi:U\times [0,T)\times
V\times[0,T)\rightarrow\mathbb{R}$ be of class $C^2.$ Assume that
${\color{black}(x,t,y,s)} \mapsto
u(x,t)-v(y,s)-{\color{black}\Phi}(x,t,y,s)$ reaches a local maximum
in $(\overline{x},\overline{t},\overline{y},\overline{s})\in U\times
[0,T)\times V\times[0,T).$ Computing the following quantities at the
point $(\overline{x},\overline{t},\overline{y},\overline{s})$, we
set $\tau_1=\partial_t\Phi$, $\tau_2=-\partial_s\Phi$,
$p_1=D_x\Phi$, $p_2=-D_y\Phi$ and
\begin{equation}\label{eq::Abar}
\overline{A}=
\left(\begin{array}{cc}
A & A_1 \\
A_1^T & A_2 \\
\end{array}\right)
\quad \mbox{with}\quad
{A}=
\left(\begin{array}{cc}
D^2_{xx}\Phi & D^2_{xy}\Phi \\
D^2_{yx}\Phi & D^2_{yy}\Phi \\
\end{array}\right),
\quad
A_1=
\left(\begin{array}{cc}
D^2_{xt}\Phi & D^2_{xs}\Phi \\
D^2_{yt}\Phi & D^2_{ys}\Phi \\
\end{array}\right),
\quad
A_2=
\left(\begin{array}{cc}
D^2_{tt}\Phi & D^2_{ts}\Phi \\
D^2_{st}\Phi & D^2_{ss}\Phi \\
\end{array}\right).
\end{equation}
Let $\overline{I}$ be the identity matrix of $\R^{2(n+1)}$. Then for
every $\mu>0$ such that $\mu \overline{A}<\overline{I},$ there
exists $X,Y\in {\mathcal S}^n$ such that (where we denote by $*$
some elements that we do not precise):
$$(\tau_1,p_1,
\left(\begin{array}{cc}
X & * \\
* & * \\
\end{array}\right)
)\in\overline{\mathcal{D}}^+\overline{u}(\overline{x},\overline{t}),
\;\;(\tau_2,p_2,
\left(\begin{array}{cc}
Y & * \\
* & * \\
\end{array}\right)
)\in\overline{\mathcal{D}}^-\overline{v}(\overline{y},\overline{s}),$$
and
$$-(\frac{1}{\mu}+\|\overline{{A}}\|)\hat{I}\leq
\left(\begin{array}{cc}
X & 0 \\
0 & -Y \\
\end{array}\right)
\leq {A}+2\mu({A}^2+A_1\cdot A_1^T).$$
where $\hat{I}$ is the identity matrix of $\R^{2n}$.
\end{corollary}

\noindent \proof $\;$ Because $\Phi\in C^2$ and
$u(x,t)-v(y,s)-\Phi(x,t,y,s)$ admits a local maximum in
$(\overline{x},\overline{t},\overline{y},\overline{s})$, we can then
apply the elliptic Ishii's Lemma (Lemma $\ref{elliptic}$) with the
new variables $\tilde{x}:=(\overline{x},\overline{t})$ and
$\tilde{y}:=(\overline{y},\overline{s})$. We obtain, for every $\mu$
satisfying $\mu \tilde{A}<\overline{I}$ with the $2(n+1)\times
2(n+1)$ matrix $\tilde{A} =D^2\Phi$, that there exists
$\tilde{X},\tilde{Y}\in {\mathcal S}^{n+1}$ with
$$\tilde{X}=
\left(\begin{array}{cc}
X & C \\
{}^tC & \delta \\
\end{array}\right),\quad
\tilde{Y}=
\left(\begin{array}{cc}
Y & D \\
{}^tD & \rho \\
\end{array}\right)
\quad \mbox{with}\quad X,Y\in \S^n, \quad
C,D\in\mathbb{R}^n
\quad \mbox{and}\quad \delta,\rho\in\mathbb{R}$$
such that:
$$\left(\left(p_1,\tau_1\right),\tilde{X}\right)
\in\overline{\mathcal{D}}^+u(\overline{x},\overline{t}),\quad
\left(\left(p_2,\tau_2\right),\tilde{Y}\right)
\in \overline{\mathcal{D}}^-v(\overline{y},\overline{s})$$
and
\begin{equation}\label{ineq}
    -(\frac{1}{\mu}+\|\tilde{A}\|)\overline{I}\leq\left(
                 \begin{array}{c|c}
                   {\begin{array}{cc}
                     X & C \\
                     {}^tC & \delta \\
                   \end{array}}&{ 0} \\
                   \hline
                   {0} &{\begin{array}{cc}
                                      -Y & -D \\
                                      -{}^tD & -\rho \\
                                    \end{array}}\\
                 \end{array}
               \right)\leq \tilde{A}+\mu \tilde{A}^2.
\end{equation}
We first remark that the matrix $\overline{A}$ (defined in
(\ref{eq::Abar})) is obtained from $\tilde{A}$ by relabeling the
vectors of the basis (going from coordinates $(x,t,y,s)$ for
$\tilde{A}$ to coordinates $(x,y,t,s)$ for $\overline{A}$).
Therefore we have $ {\color{black}\|}\tilde{A} {\color{black}\|}=
{\color{black}\|}\overline{A} {\color{black}\|}$ and the condition
$\mu\tilde{A} < \overline{I}$ is equivalent to
$\mu  \overline{A}< \overline{I}$.\\
Next, for
$\xi,\eta\in\mathbb{R}^n,$ applying the vector $\tilde{V}=(\xi,0,\eta,0)^T$ to
the matrix inequality $(\ref{ineq}),$ yields with $V=(\xi,\eta)^T$:
\begin{eqnarray}\label{eq:9.4}
\displaystyle{-(\frac{1}{\mu}+\|\overline{A}\|)
<\hat{I} \cdot V,V>}&\leq&
\displaystyle{
<\left(\begin{array}{cc}
X & 0 \\
0 & -Y \\
\end{array}\right)
\cdot V,V>}\nonumber\\
&\leq&<(\tilde{A}+\mu
   \tilde{A}^2)\cdot \tilde{V},\tilde{V}> {\color{black}=:} J.
\end{eqnarray}
Remark that using the relabeling of the vectors of the basis, we have
$$<\tilde{A}^k\tilde{V},\tilde{V}>
= <\overline{A}^k \overline{V},\overline{V}>
\quad\mbox{with}\quad
\overline{V}=\left(\begin{array}{l}
V\\
0_{\R^2}
\end{array}\right),\quad \mbox{for}\quad k=1,2$$
We compute
$$\left(\begin{array}{l}
V\\
0
\end{array}\right)^T
\left(\begin{array}{cc}
A & A_1\\
A_1^T & A_2
\end{array}\right)
\left(\begin{array}{l}
V\\
0
\end{array}\right)
= V^T A V$$
and
$$\left(\begin{array}{l}
V\\
0
\end{array}\right)^T
\left(\begin{array}{cc}
A & A_1\\
A_1^T & A_2
\end{array}\right)^2
\left(\begin{array}{l}
V\\
0
\end{array}\right)
= V^T A^2 V+  V^T A_1 A_1^T V + 2 < A V, A_1^T V> \le 2\left(|AV|^2 + |A_1^T V|^2\right)$$
This gives
$$J\le  V^T A V + 2\mu \left(|AV|^2 + |A_1^T V|^2\right)$$
and then (\ref{eq:9.4}) implies
$$-(\frac{1}{\mu}+\|\overline{A}\|)\hat{I}\leq
\left(\begin{array}{cc}
X & 0 \\
0 & -Y \\
\end{array}\right)
\leq A+2\mu(A^2+A_1\cdot A_1^T),$$
which achieves the proof. $\hfill{\Box}$

\bigskip
\noindent {\bf Acknowledgments}\\
{\color{black} The first two authors highly appreciate the kind
hospitality of Ecole des Ponts ParisTech (ENPC) while starting the
preparation of this work. The third author would like to thank M.
Jazar for welcoming him as invited professor at the LaMA-Liban for
several periods during the preparation of this work. Part of this
work was supported by the ANR MICA project $(2006-2010).$}


\end{document}